\documentclass[11pt,reqno]{amsart}
\usepackage[T1]{fontenc}
\usepackage[utf8]{inputenc}
\usepackage{lmodern}
\usepackage[a4paper,margin=27mm]{geometry}
\usepackage{amsmath,amssymb,amsthm,mathtools}
\usepackage{microtype}
\usepackage{booktabs,array,enumitem}
\usepackage{tikz}
\usetikzlibrary{arrows.meta,calc,patterns}
\usepackage{xcolor}
\usepackage{hyperref}
\hypersetup{colorlinks=true,linkcolor=blue!40!black,citecolor=blue!40!black,urlcolor=blue!40!black,pdftitle={Division in quasianalytic classes and sharp growth thresholds}}
\numberwithin{equation}{section}
\newtheorem{theorem}{Theorem}[section]
\newtheorem{proposition}[theorem]{Proposition}
\newtheorem{lemma}[theorem]{Lemma}
\newtheorem{corollary}[theorem]{Corollary}
\theoremstyle{definition}
\newtheorem{definition}[theorem]{Definition}
\newtheorem{remark}[theorem]{Remark}

\newcommand{\R}{\mathbb R}
\newcommand{\C}{\mathbb C}
\newcommand{\N}{\mathbb N}
\newcommand{\CM}{C^M}

\DeclareMathOperator{\dist}{dist}
\DeclareMathOperator{\ord}{ord}

\setlist[enumerate]{label=\textup{(\roman*)},leftmargin=*,itemsep=3pt,topsep=5pt}
\allowdisplaybreaks[2]
\title[Division in quasianalytic DC classes]{Division and noetherianity in quasianalytic Denjoy-Carleman classes}

\author[A.~Belotto da Silva]{Andr\'e Belotto da Silva}
\author[A.~Kiro]{Avner Kiro}
\address[A.~Belotto da Silva]{Universit\'e Paris Cit\'e and Sorbonne Universit\'e, UFR de Math\'ematiques, Institut de Math\'ematiques 
de Jussieu-Paris Rive Gauche, UMR7586, F-75013 Paris, France, and Institut universitaire de France (IUF)}
\email{andre.belotto@imj-prg.fr}
\address[A.~Kiro]{Department of Mathematics,
Weizmann Institute of Science,
Rehovot 7610001, Israel}
\email{avner-ephraiem.kiro@weizmann.ac.il}
\thanks{Research supported by Agence Nationale de la Recherche (ANR) projects ANR-25-ERCC0003-01 (Belotto da Silva) and by the European Union (ERC, SharpOS, 101087910), and by the Israel Science Foundation (grant No. 2067/23) (Kiro).}
\date{\today}

\begin{document}
\begin{abstract}
We show that the local rings of quasianalytic Denjoy–Carleman classes need not be Noetherian. More precisely, we construct quasianalytic Denjoy–Carleman classes $\CM$ whose local rings $\CM(\mathbb{R}^d,0)$ are not Noetherian for any $d\geq 2$. Conversely, we exhibit quasianalytic classes whose local rings in dimension two are regular and Noetherian.

Our approach is based on the problem of exact division: if $f=\varphi g$, where $f\in\CM$, $\varphi$ is real-analytic, and $g\in C^\infty$, under which conditions does it follow that $g\in\CM$? We introduce a quantitative scale-growth condition on the defining weight sequence $M$ and prove that exact division holds in two variables whenever this condition is sufficiently strong relative to the rate at which the nonreal branches of the complex zero set of $\varphi$ approach the real plane. We also construct quasianalytic weight sequences with prescribed growth thresholds and establish corresponding counterexamples to exact division. These yield nonclosed principal ideals and, consequently, non-Noetherian local rings.
\end{abstract}
\maketitle

\section{Introduction}\label{sec:intro}

Quasianalytic Denjoy--Carleman classes $\CM$ --parametrised by the weight sequence $M = (M_k)_{k\in \mathbb{N}}$ satisfying certain natural qualitative assumptions, see Definition \ref{def:QDC} -- form an important family of
$C^\infty$ functions, with origins in classical analysis
\cite{Borel1900,Denjoy1921,Carleman1926,Mandelbrojt1952}, and with
important connections to analysis and o-minimal geometry
\cite{RolinSpeisseggerWilkie2003}\footnote{For some developments and applications in analysis, see, for
instance,
\cite{Dynkin1980,ChaumatChollet1999,Thilliez2001,
Thilliez2002Sharp,FurdosNenningRainerSchindl2020,Rainer2022Survey}; for connections with o-minimal geometry and smooth definable
quasianalytic functions, see
\cite{Miller1995,RolinSpeisseggerWilkie2003,Nowak2013,
BelottoBiborskiBierstone2017,BelottoBierstoneChow2018,
BelottoBierstone2023}.}.
Although they share many of the basic local properties of real-analytic
functions---the Taylor homomorphism is injective, they are closed under
differentiation and composition, and admit inverse and implicit function
theorems; see Remark \ref{rem:QDCbasic-properties}---their local algebra is much
less understood.

A basic longstanding question, already raised by Childress \cite{Childress1976} and explicitly formulated by Bierstone and Milman in their study of quasianalytic geometry \cite[\S6]{BierstoneMilman2004}, is whether the local rings of germs in quasianalytic Denjoy--Carleman classes are Noetherian; see also \cite[\S3.6]{Thilliez2008}. Edward Bierstone brought this problem to our attention, and our discussions and earlier attempts to address it together have been an important motivation for the present work. Noetherianity is one of the fundamental finiteness properties of commutative algebra and underlies
much of classical local algebraic and analytic geometry.  For analytic
germs it follows, in particular, from Weierstrass preparation and
division.  The quasianalytic problem therefore asks to what extent this
basic algebraic feature of analytic functions persists when analyticity is replaced by quasianalyticity.

The classical Weierstrass machinery does not extend to proper quasianalytic Denjoy--Carleman classes.  Childress \cite{Childress1976} showed that Weierstrass division fails in this setting, while
Acquistapace, Broglia, Bronshtein, Nicoara and Zobin
\cite{AcquistapaceEtAl2014} proved failure of the
Weierstrass preparation theorem.  More generally, Parusi\'nski and Rolin \cite{ParusinskiRolin2014} showed that, for a quasianalytic system satisfying
the usual stability properties, Weierstrass preparation can hold only in
the analytic case.  These results rule out the classical Weierstrass
route to Noetherianity, but do not by themselves determine whether
quasianalytic Denjoy--Carleman local rings can be Noetherian. Our first result shows that there is in fact no uniform answer:

\begin{theorem}\label{thm:intro-noetherianity}
There exists a quasianalytic Denjoy--Carleman class $\CM$, with $M$ of moderate growth, such that the local ring of germs $\CM(\mathbb{R}^d,0)$ is not Noetherian for every $d\geq 2$.

On the other hand, for each fixed $0<a\leq 1$, the weight
\[
    M_k=(\log(k+e))^{ak},
\]
defines a class whose local ring $\CM(\mathbb{R}^2,0)$ is a regular Noetherian local ring of dimension two.
\end{theorem}

In particular, the standard qualitative assumptions on a quasianalytic Denjoy--Carleman class do not by themselves determine the algebraic behavior of its local ring. Finer quantitative properties of the defining weight sequence $M$ play an essential role. In fact, both assertions in Theorem~\ref{thm:intro-noetherianity} are consequences of the main investigation, and technical heart, of the paper: exact division without loss of regularity in quasianalytic Denjoy–Carleman classes.

\begin{remark}\label{rem:novelty_buildingblock}
Beyond the counterexamples to Noetherianity, our proof introduces a construction method of independent interest. A classical approach to producing functions with prescribed derivative behavior is to sum carefully chosen building blocks, often involving rational functions. This approach has a long history, going back to ideas of Borel \cite{Borel1900}; see also \cite{AcquistapaceBrogliaNicoara2016,Thilliez2002Sharp,BelottoBierstoneKiro2022}. Adapting such constructions to genuinely multidimensional phenomena is difficult, and in our earlier work with Bierstone \cite{BelottoBierstoneKiro2022}, we obtained examples exhibiting intrinsically two-dimensional behavior.

Here we introduce a different type of building block: polynomial approximations of smooth localized quotients constructed using cutoff functions, with quantitative estimates compatible with quasianalytic regularity. The underlying polynomial-replacement mechanism is not intrinsically two-dimensional and may provide a useful method for constructing examples in higher-dimensional quasianalytic analysis.
\end{remark}

\subsection{Exact division}

The main problem addressed in this paper is the following. Let $\varphi$ be a
nonzero real-analytic germ and suppose that
\[
        f=\varphi g, \qquad f\in \CM, \quad g\in C^\infty .
\]
Under which conditions does it follow that $g$ belongs to the same class
$\CM$?

The regularity of such quotients has been extensively studied in Denjoy–Carleman classes. Chaumat and Chollet
\cite{ChaumatChollet1999,ChaumatChollet2004} obtained several division
results in this setting; in particular, their hyperbolic division theorem
shows that division by a hyperbolic polynomial preserves the Denjoy-Carleman class, in both the quasianalytic and nonquasianalytic cases. For general analytic plane divisors, a particularly precise
answer in two variables was obtained by Thilliez
\cite{Thilliez2001,Thilliez2002Sharp,Thilliez2002Formal}. For a strongly
regular, \emph{nonquasianalytic} Denjoy--Carleman class and a real-analytic
plane germ $\varphi$, he showed that the loss of regularity under division is determined by the rate at which the nonreal complex branches of $\{\varphi=0\}$ approach the real plane, that is,
\[
       g\in C^{M^{d(\varphi)}},
       \qquad
       M^{d(\varphi)}=(M_k^{\,d(\varphi)})_{k\geq0},
\]
where $d(\varphi)\geq1$ is a {\L}ojasiewicz-type separation exponent that can
be read from the Puiseux expansions of the complex branches. This idea plays a central role in this work as well.

Here we work in the quasianalytic setting and ask when there is \emph{no loss of regularity}. Our main results give a sharp answer to the exact division problem in dimension two, see Theorems \ref{thm:positive} and \ref{thm:sharpnessIntro} below. General quasianalytic division and solvability results, based on resolution of singularities, show that one can obtain a quotient in a larger quasianalytic class, with a controlled loss of regularity; see
\cite{BelottoBiborskiBierstone2017,BierstoneMilman2017}.
These geometric methods, however, intrinsically lead to such a loss of regularity \cite{BelottoBierstoneKiro2022}. Determining when the quotient remains in the original class therefore requires a finer control at the level of Denjoy-Carleman estimates.

\subsubsection{Quasianalytic Denjoy-Carleman classes} In order to state our main results, we now recall the basic notions of Denjoy-Carleman classes.

\begin{definition}[Quasianalytic Denjoy-Carleman classes]\label{def:QDC}
Let $M=(M_k)_{k\geq0}$ be an increasing logarithmically convex sequence of positive numbers, normalized by $M_0=1$; that is,
\[
1\leq M_k\leq M_{k+1},
\qquad
M_k^2\leq M_{k-1}M_{k+1}\quad(k\geq1).
\]
We assume throughout that
\begin{equation}\label{eq:standing-diff}
    M_{k+1}\leq D^{k+1}M_k \qquad (k\geq0)
\end{equation}
for some $D\geq1$, and that
\begin{equation}\label{eq:standing-qa}
    \sum_{k\geq0}\frac{M_k}{(k+1)M_{k+1}}=\infty .
\end{equation}

For an open set $U\subset\R^d$, the Denjoy--Carleman class $\CM(U)$ consists of the functions $f\in C^\infty(U)$ such that, for every compact $K\Subset U$, there are constants $A,B>0$ satisfying
\begin{equation}\label{eq:class-bound}
    \sup_{x\in K}|\partial^\alpha f(x)|
    \leq AB^{|\alpha|}|\alpha|!M_{|\alpha|}
    \qquad(\alpha\in\N^d).
\end{equation}
\end{definition}

\begin{remark}[Basic properties of quasianalytic DC classes]\label{rem:QDCbasic-properties}
The condition \eqref{eq:standing-qa} is the Denjoy--Carleman quasianalyticity criterion, while \eqref{eq:standing-diff} implies that $\CM$ is closed under differentiation. Under these hypotheses, $\CM$ is also closed under composition, division by a coordinate, and the inverse and implicit function theorems; see, for example, \cite[\S\S2--4]{BierstoneMilman2004}, and the references therein for the classical sources. In particular, a germ in $\CM$ is uniquely determined by its Taylor series.

If $\sup_k M_k^{1/k}<\infty$, then $\CM$ coincides with the class of real-analytic functions, where the division problem considered here has a positive answer by analytic division. In the sequel, unless explicitly stated otherwise, we therefore restrict to the genuinely nonanalytic case,
\begin{equation}\label{eq:proper-weight}
    M_k^{1/k}\longrightarrow\infty .
\end{equation}
\end{remark}
\begin{remark}[Convention for the weight sequence]
We use the convention in which the factorial is not included in the weight sequence $M$; thus some references instead work with $m_k=k!M_k$. Note, moreover, that replacing $|\alpha|!$ by $\alpha!$ in \eqref{eq:class-bound} does not change the class.
\end{remark}

\subsubsection{Scale-growth condition} The first ingredient in our division theorem concerns the Denjoy-Carleman class $\CM$. In the sequel, we consider two standard associated functions of the weight sequence $M$, see e.g. \cite{Thilliez2002Sharp,Rainer2022Survey}, which provide  a convenient way to encode its quantitative growth:
\begin{equation}\label{eq:associated}
 h_M(t)=\inf_{k\geq0} t^k M_k,\qquad \text{ and }  \qquad \omega_M(T)
    =-\log h_M(1/T)
    =\sup_{k\geq0}\log\frac{T^k}{M_k},
\end{equation}
for $t,T>0$; the last function is used only for sufficiently large $T$. Moreover, motivated by the role of the associated function $h_M$ in the sharp division estimates of Thilliez \cite{Thilliez2002Sharp}, we consider
\begin{equation}\label{eq:associated2}
\Lambda_M(T)=\log \omega_M(T)
    =\log\bigl(-\log h_M(1/T)\bigr),
\end{equation}
for all sufficiently large $T$. This last function is used to formulate the growth condition that controls exact division:

\begin{definition}[Scale-growth condition]\label{def:Gbeta}
Let $0\leq\beta<1$. We say that $M$ satisfies the scale-growth condition $(G_\beta)$ if there exist $R>1$ and $c>0$ such that
\begin{equation}\label{eq:Gbeta}
    \Lambda_M(T)-\Lambda_M(T/R)\geq cT^\beta
    \qquad (T\gg1).
\end{equation}
\end{definition}

\begin{remark}[Basic properties of the scale-growth condition]\label{rem:Gbeta-basic}
The scale-growth conditions $(G_\beta)$ have the following basic properties.
\begin{enumerate}
\item Condition $(G_\beta)$ can equivalently be written as
\[
    h_M(t)^{\exp(-ct^{-\beta})}\leq h_M(Rt)
\]
for all sufficiently small $t>0$, after possibly changing the constants.

\item Condition $(G_0)$ is equivalent to moderate growth of $M$, that is, there exists $H\geq1$ such that
\[
    M_{j+k}\leq H^{j+k}M_jM_k
    \qquad(j,k\geq0).
\]

\item If $M$ satisfies $(G_\beta)$, then it satisfies $(G_{\beta'})$ for every
$0\leq\beta'\leq\beta$.

\item The conditions $(G_\beta)$ depend only on the Denjoy--Carleman class and not on the particular choice of an equivalent weight sequence. More precisely, they are invariant under exponential equivalence: if
\[
    C^{-k}M_k\leq N_k\leq C^kM_k
    \qquad(k\geq0)
\]
for some $C\geq1$, then $M$ satisfies $(G_\beta)$ if and only if $N$ does.
\end{enumerate}
Assertions (2)--(4) are proved in Lemma~\ref{lem:weight-calculus}.
\end{remark}

\subsubsection{The contact exponent $q$} The second ingredient in our division theorem is geometric. Following Thilliez, we associate with a nonzero real-analytic plane germ $\varphi$ vanishing at the origin an exponent $q(\varphi)$ measuring the rate at which the nonreal branches of its complex zero set approach the real plane. More explicitly, after passing to the reduced divisor and choosing adapted real coordinates, write the Puiseux roots over the two real half-axes as
\(
    y=\lambda_\nu(\pm t),
\)
where $\nu=1,\ldots,r$ for some $r\in \mathbb{N}$. For each nonreal branch, there exists $q_{\nu}\in \mathbb{Q}$ such that $|\operatorname{Im}\lambda_\nu(\pm t)|\asymp t^{q_\nu}$, and we set:
\[
q(\varphi)=\max\bigl(\{1\}\cup\{q_\nu\}\bigr),
\]
where the maximum is taken in respect to all imaginary roots. Moreover, we set $q(\varphi)=1$ when $\varphi$ is an analytic unit. This is the exponent denoted $d(\varphi)$ by Thilliez; see \cite[\S\S1.3--1.4]{Thilliez2002Sharp}. Section~\ref{sec:geometry} discusses in more detail its geometric interpretation and properties.

\subsubsection{The division Theorem} We are finally ready to state the first main Theorem:

\begin{theorem}[Exact division in dimension two]\label{thm:positive}
Let $0\leq\beta<1$, and assume that $M$ satisfies $(G_\beta)$. Let $\varphi$ be a nonzero real-analytic germ in two variables such that
\begin{equation}\label{eq:geometric-threshold}
    q(\varphi)\leq \frac{1}{1-\beta}.
\end{equation}
If two germs $f\in \CM_0(\R^2)$ and $g\in C^\infty_0(\R^2)$ satisfy
\(
    f=\varphi g,
\)
then $g\in \CM$.
\end{theorem}

Before discussing the proof, we record an immediate corollary on formal division, obtained by combining Theorem~\ref{thm:positive} with \cite[Proposition~4.9]{BelottoBiborskiBierstone2017}, which is particularly useful in applications:

\begin{corollary}[Formal division]\label{cor:formal}
Under the hypotheses of Theorem~\ref{thm:positive}, let $f\in \CM_0(\R^2)$. If the Taylor series $\widehat f$ is formally divisible by $\widehat\varphi$, say
\[
    \widehat f=\widehat\varphi\,G
\]
for some $G\in\R[[x,y]]$, then there exists a unique
$g\in \CM_0(\R^2)$ such that
\[
    f=\varphi g,
    \qquad
    \widehat g= G.
\]
\end{corollary}
\begin{proof}
It follows from \cite[Proposition~4.9]{BelottoBiborskiBierstone2017} that there exists a smooth function $g$ such that $f=\varphi g$ and $\hat{g} =G$. We conclude from Theorem~\ref{thm:positive}.
\end{proof}

The proof of Theorem~\ref{thm:positive}, given in Section~\ref{sec:positive}, has two main new quantitative ingredients: Lemma~\ref{pos:branch-smallness} and Proposition~\ref{pos:local-division}. The elementary Lemma~\ref{pos:small-product} also plays an essential role in the argument. We briefly explain their respective roles and how they relate to the classical methods on which our proof builds.

The proof relies on holomorphic approximation methods for Denjoy-Carleman  functions, originating in Dynkin's theory of almost analytic extensions \cite{Dynkin1980,Dynkin1993}. In this approach, a function $u\in\CM$ is approximated by holomorphic functions $u_s$, defined on complex neighbourhoods of size comparable to $s>0$, with approximation errors bounded by $Ch_M(As)$, where $A,C>0$ are constants independent of $s$. Conversely, the existence of such approximations, with suitable uniform estimates, characterizes Denjoy--Carleman regularity. Propositions~\ref{prop:hol-approx} and~\ref{prop:approx-characterization} are variants of standard approximation results in this theory, whose proofs are included in Appendix~\ref{app:holomorphic}.

The main difficulty is to carry out division while preserving these quantitative estimates. After reducing to a Weierstrass polynomial $W$, write $f=Wg$. One of the central ideas is to work at two different scales: we seek holomorphic approximations to the quotient $g$ at scale $\sigma$, starting from an approximation $F_s$ of the numerator at the larger scale $s=L\sigma$, where $L>1$ is a sufficiently large fixed constant. Although $f$ is divisible by $W$, its holomorphic approximant $F_s$ need not vanish on the complex zero set $\{W=0\}$. Our strategy is to control this defect of divisibility by obtaining bounds of the form $Ch_M(A\sigma)$ for the values of $F_s$ on the relevant parts of the complex divisor.

We first establish such estimates near the origin (Lemma~\ref{pos:central}), using known bounds for formal quotients \cite{Thilliez2002Formal} and optimal Taylor truncation. The next step is to propagate these estimates along the complex branches of the divisor. Their geometry is described by the separation exponent studied by Thilliez \cite{Thilliez2002Sharp}, which determines how far the relevant branches may extend while remaining close to the real plane. The difficulty is that analytic continuation weakens the smallness estimates, and the number of continuation steps increases as $s$ tends to zero. The key contribution of Lemma~\ref{pos:branch-smallness} is to control this accumulated loss while retaining an estimate of the form $Ch_M(A\sigma)$. Lemma~\ref{pos:absorption} shows how the scale-growth condition $(G_\beta)$ compensates for the loss through the choice of the enlargement factor $L$. This two-scale mechanism is essential for obtaining the result also at the critical exponent $q(\varphi)=(1-\beta)^{-1}$.

With these estimates on the complex divisor, Proposition~\ref{pos:local-division} constructs local holomorphic approximate quotients, with polynomial bounds on their complex norms and division remainders bounded by $Ch_M(A\sigma)$. The construction relies on classical Cauchy division, but requires a careful choice of integration contours and quantitative control of the remainders, since $F_s$ is not exactly divisible by $W$. Finally, Lemma~\ref{pos:small-product} provides an elementary estimate allowing us to pass from the smallness of $W(g-G)$ to that of $g-G$, even near the real zero set of $W$. The resulting bounds, after absorbing fixed losses, give holomorphic approximations to the smooth quotient $g$ with the required derivative control. Proposition~\ref{prop:approx-characterization} then yields $g\in\CM$.

\subsubsection{Sharpness of exact division} The next result shows that the geometric threshold in Theorem~\ref{thm:positive} is sharp, by existence of counterexamples (it does not assert a converse for an arbitrary fixed weight failing $(G_\beta)$):

\begin{theorem}[Sharpness]\label{thm:sharpnessIntro}
For every $0<\beta_0<1$, there exists a quasianalytic Denjoy--Carleman class $\CM$ such that
\begin{equation}\label{eq:exact-spectrum}
    \{\beta\in[0,1):(G_\beta)\text{ holds}\}=[0,\beta_0).
\end{equation}
Moreover, for every pair of integers $1\leq m<n$ satisfying
\(
    1-\frac{m}{n}\geq\beta_0,
\)
there exist germs
\[
    f\in \CM_0(\R^2),
    \qquad
    g\in C^\infty_0(\R^2)\setminus \CM_0(\R^2),
\]
such that
\(
    f=(y^{2m}+x^{2n})g.
\)
\end{theorem}

A more precise quantitative version of Theorem~\ref{thm:sharpnessIntro} is given in Theorem~\ref{thm:sharpness} below. Besides an explicit growth estimate for the weight sequence, it establishes that the numerators of the counterexamples can be approximated by polynomial multiples of the corresponding divisor, with convergence in a fixed Denjoy--Carleman norm.

The weight sequence $M$ is constructed in Section~\ref{sec:staircase}, and the counterexamples in Section~\ref{sec:negative}. As discussed in Remark \ref{rem:novelty_buildingblock}, our proof follows a classical building-block strategy. More specifically, we organize the construction along the same lines as in our earlier work with Bierstone \cite{BelottoBierstoneKiro2022}, with two main new ingredients: a new type of building blocks, based on quantitative polynomial approximation, and a staircase weight sequence with a prescribed sharp growth threshold.

Our starting point is a simple construction using cutoff functions. This approach is natural in nonquasianalytic classes, particularly Gevrey classes, where localization by cutoff functions is available. For divisors of the form $P=y^{2m}+x^{2n}$, we consider smooth functions obtained by dividing a localized cutoff by $P$ (Lemma \ref{neg:cutoffs}). Their products with $P$ are particularly simple, while their transverse derivatives can be very large, reflecting the proximity of the complex zeros of $P$ to the real plane. However, a direct summation of these localized blocks cannot yield a nonzero quasianalytic numerator, since the resulting function would be flat at the origin.

Our main new idea is therefore to use these smooth functions only as intermediate models, replacing each of them by a polynomial before taking the sum. This is achieved by Gaussian smoothing followed by Taylor truncation (Lemma~\ref{neg:gaussian-approximation}). Although these approximation procedures are classical, the difficulty is to control simultaneously the approximation error in a Denjoy--Carleman norm, the polynomial degree, and the large derivatives responsible for the failure of division. The key observation is that the approximations can be carried out using auxiliary nonquasianalytic weights that agree with the original weight $M$ up to a prescribed order (Lemma \ref{neg:temporary-properties}). By keeping the degrees of the polynomial numerators within this range, we transfer the required estimates to the original quasianalytic class. These ingredients lead to the polynomial building blocks of Proposition \ref{neg:polynomial-blocks}.

A further ingredient is the construction, in Section~\ref{sec:staircase}, of a staircase weight sequence with a prescribed growth threshold. The idea of using long intervals of controlled growth to construct quasianalytic weights already appeared in our earlier work with Bierstone \cite[\S4]{BelottoBierstoneKiro2022}, where a simpler construction was used to separate Denjoy--Carleman classes. Here we develop this idea further, alternating growth intervals and plateaux in the quotient sequence $M_k/M_{k-1}$. The plateaux ensure quasianalyticity (Proposition \ref{stair:regularity}), while the construction realizes a prescribed threshold for the scale-growth condition $(G_\beta)$ (Proposition \ref{stair:threshold}). At the same time, the growth intervals must be sufficiently long to accommodate the polynomial approximations described above (Lemma \ref{stair:separation}). The quantitative choices are particularly delicate at the critical exponent, where the leading exponential gain disappears and additional logarithmic estimates are needed to ensure that the normalized derivatives remain unbounded (Lemma  \ref{neg:small-estimates} and Remark \ref{neg:critical-remark}).

\subsection{Closed ideals and Noetherianity} 
As we already anticipated, the problem of exact division is closely related to the question of Noetherianity of quasianalytic local rings. As discussed by Thilliez \cite[\S3.6]{Thilliez2008}, a classical approach to proving non-Noetherianity is to exhibit an ideal which is \emph{not} closed in the local ring $\CM(\mathbb{R}^d,0)$. 

To make this precise, for a weight sequence $M$, let
\[
A_{M,d}:=\CM_0(\mathbb R^d),\qquad
B_d:=\mathbb R[[x_1,\ldots,x_d]] 
\]
denote respectively the ring of Denjoy–Carleman germs at the origin and the ring of formal power series, which is the completion of $A_{M,d}$. 

\begin{definition}\label{def:closedideal}
Denote by \(\mathfrak m\) the maximal ideal of \(A_{M,d}\). An ideal \(I\subset A_{M,d}\) is said to be \emph{closed} if
\[
I=\bigcap_{k\geq 1}(I+\mathfrak m^k).
\]
\end{definition}

\begin{remark}\label{rem:ClosedIdealsDef}
When \(M\) is quasianalytic, the Taylor expansion map $T : A_{M,d} \to B_d$ is injective, and the image is identified as a subring of $B_d$. Given an ideal $I \subset A_{M,d}$, denote by $I B_d$ the ideal of $B_d$ generated by all elements in $T(I) \subset B_d$. Then $I$ is closed if, and only if:
\[
I=IB_d\cap A_{M,d}.
\]
Thus, an ideal is closed if formal membership implies membership in the quasianalytic ring.
\end{remark}

Theorem~\ref{thm:sharpnessIntro} also has consequences
for the closedness of ideals. Indeed, the counterexamples
to exact division yield non-closed principal ideals.
Since every ideal in a Noetherian local ring is closed,
by Krull's intersection theorem, we obtain:

\begin{corollary}\label{cor:nonnoetherian}
For every weight $M$ in Theorem~\ref{thm:sharpnessIntro}
and every pair of integers $1\leq m<n$ satisfying
\(
1-\frac{m}{n}\geq\beta_0,
\)
the principal ideal generated by
$P=y^{2m}+x^{2n}$ in $A_{M,2}$ is not closed.
Consequently, $A_{M,d}$ is not Noetherian
for any $d\geq2$.
\end{corollary}
\begin{proof}
Let $f=Pg$ be as in Theorem~\ref{thm:sharpnessIntro}.
Then $\widehat f\in \widehat P B_2$, but $f\notin PA_{M,2}$,
since otherwise $g$ would belong to $A_{M,2}$.
Thus $PA_{M,2}$ is not closed. By Krull's intersection theorem, $A_{M,2}$ is not Noetherian. For $d\geq2$, restriction to a coordinate plane gives a surjective homomorphism $A_{M,d}\to A_{M,2}$, which proves the remaining assertion.
\end{proof}

This result implies the negative statement of Theorem \ref{thm:intro-noetherianity}.

\begin{remark}
With $\beta_0=1/2$, the hypothesis of Theorem~\ref{thm:sharpnessIntro} allows one to consider $P=y^2+x^4$. It therefore gives an explicit nonclosed principal ideal of the form considered in \cite[Remark~4.2.2]{Thilliez2008}.
\end{remark}

On the positive side, Corollary~\ref{cor:formal} shows
that, if the weight sequence $M$ satisfies the condition $(G_{\beta})$ for all $0\leq \beta <1$, then principal ideals generated by real-analytic germs are closed.
In dimension two, we can obtain a stronger result: under the same assumption, we can show that every principal ideal of $A_{M,2}$ is closed, see Proposition \ref{prop:all-principal-contraction}. Combining this with the well-known local-algebraic argument of \cite[Proposition~3]{Valabrega1973}, we obtain:

\begin{theorem}\label{thm:plane}
Suppose that $M$ is a weight sequence that satisfies $(G_\beta)$ for every
$0\leq\beta<1$. Then $A_{M,2}$ is a regular
(in particular, Noetherian) local ring of dimension two.
\end{theorem}
\begin{remark}
Under the same hypothesis of Theorem \ref{thm:plane}, it follows that every ideal of $A_{M,2}$ is closed.
\end{remark}

The proof of Theorem~\ref{thm:plane} is in Section~\ref{sec:algebra}. Technically, an important ingredient in the proof is a finite-determinacy Lemma \ref{lem:reduced-normal-form} for formally reduced plane germs: every such germ is equivalent, under a local Denjoy--Carleman change of coordinates, to a sufficiently high Taylor polynomial. This is a Denjoy--Carleman version of a classical finite-determinacy theorem, see e.g. \cite[Chapter VIII, §5, Proposition 5.1.2]{Tougeron1972}, obtained using the stability of Denjoy--Carleman classes under flows \cite{RainerSchindl2016}. In particular, it does not require the scale-growth condition $(G_\beta)$.

Finally, Proposition~\ref{prop:log-weights} shows that the weight sequence:
\[
 M_k=(\log(k+e))^{ak},\qquad 0<a\leq1;
\]
satisfies the property $(G_\beta)$ for all $0\leq \beta <1$ and, therefore, gives an explicit example of a class such that $A_{M,2}$ is Noetherian, which implies the positive statement of Theorem \ref{thm:intro-noetherianity}.

\begin{remark}[Silva topology]\label{rem:Silva}
The notion of closedness in Definition~\ref{def:closedideal}
corresponds to the $\mathfrak m$-adic topology of the local
ring $A_{M,d}$. There is also a natural functional-analytic
topology on this ring, known as the \emph{Silva topology}, see details in e.g. \cite[\S3.5]{Thilliez2008}.

Since Taylor-jet maps are continuous in the Silva topology, every ideal that is closed in the $\mathfrak m$-adic topology
is also closed in the Silva topology.
Consequently, our positive closedness results hold in both topologies.

Moreover, the counterexamples of
Theorem~\ref{thm:sharpnessIntro} provide principal ideals that are not closed in either topology. We show this result in Corollary \ref{cor:closednessBoth} below.
\end{remark}

\subsection{Relation to earlier work}

Theorem~\ref{thm:positive} concerns exact division: it controls the
regularity of a smooth quotient already known to exist. Full Weierstrass
division requires a quotient and a remainder in the class for arbitrary
numerators. In a proper quasianalytic Denjoy--Carleman class, Childress
\cite{Childress1976} shows that hyperbolicity is necessary for
Weierstrass division by a distinguished polynomial; Chaumat--Chollet
\cite{ChaumatChollet2004} establish its sufficiency under the standing
weight assumptions. More generally, under the usual stability
assumptions, Elkhadiri--Sfouli \cite{ElkhadiriSfouli2008} show that a
quasianalytic system satisfying full Weierstrass division must be
analytic. The corresponding obstruction to Weierstrass preparation is
established in \cite{ParusinskiRolin2014,AcquistapaceEtAl2014}.
The example $x^2+y^2$ already illustrates the distinction between
closedness and full division: the ideal $(x^2+y^2)A_{M,2}$ is closed
in the fixed-neighborhood, functional-analytic sense of
\cite[\S4.1 and Example~4.2.1]{Thilliez2008}, although full
Weierstrass division by $x^2+y^2$ fails.

Elkhadiri \cite[Definitions~2.2 and 4.7, Corollary~5.7]{Elkhadiri2012}
proves full Weierstrass division for a well-behaved Noetherian
differentiable system. Noetherianity is required in every dimension,
and well-behavedness imposes strong injectivity for elementary
blow-up and ramification substitutions: a formal series is the Taylor
series of a germ in the system whenever its transformed series is.
The plane Noetherianity conclusion of Theorem~\ref{thm:plane}
supplies neither Noetherianity in every dimension nor these descent
properties.

In the strongly regular nonquasianalytic setting, Thilliez
\cite{Thilliez2002Sharp} obtains sharp quotient estimates in terms of
the same branchwise exponent, denoted by $d(P)$ in his notation.
Earlier and related work on closed ideals and controlled formal
division includes
\cite{Thilliez2001,Thilliez2002Formal,ChaumatChollet1997}.
Bounds for formal coefficients at a point, or at each point with
point-dependent constants, do not by themselves give the derivative
bounds uniform on a neighborhood that are required for membership
in a fixed Denjoy--Carleman class. Division and composition by
desingularization were studied by Chaumat--Chollet
\cite{ChaumatChollet1999}.

Bierstone--Milman \cite[Theorem~6.4]{BierstoneMilman2004} prove a
smooth division theorem from formal divisibility at every point.
Their later module division theorem
\cite[Theorem~1.1]{BierstoneMilman2017} gives smooth solutions of
linear equations with quasianalytic coefficients under the analogous
pointwise formal solvability condition. For a formally divisible
pair of quasianalytic germs, Belotto da Silva--Biborski--Bierstone
\cite[Proposition~4.9]{BelottoBiborskiBierstone2017} obtain a quotient
in a possibly larger quasianalytic class. Regularity loss under
composition and blowing down is treated in
\cite{BelottoBierstoneChow2018,BelottoBierstoneKiro2022}.
The present sufficient condition compensates for the quantitative
continuation loss while retaining one fixed weight.

There is also a distinction between a fixed Denjoy--Carleman class
and its smooth definable enlargement. Let $\R_M$ be the expansion
of the real field by restricted $\CM$ functions, and let
\begin{equation}\label{eq:definable-ring}
 \mathcal D_d(M)=
 \{u\in C^\infty_0(\R^d):
   u\text{ has a representative definable in }\R_M\}.
\end{equation}
Restricted functions are restrictions to compact boxes of functions
of class $\CM$ on neighborhoods of those boxes, extended by zero
outside. The structure $\R_M$ is polynomially bounded and o-minimal
\cite{RolinSpeisseggerWilkie2003}; its smooth definable germs are
quasianalytic \cite{Miller1995}. We have
$A_{M,d}\subseteq\mathcal D_d(M)$.

For $f,p\in\mathcal D_d(M)$ with $p\ne0$, Nowak
\cite{Nowak2013} proves that formal divisibility at the origin
implies divisibility in the smooth definable ring:
\[
 \widehat f\in\widehat p\,\R[[x_1,\ldots,x_d]]
 \quad\Longrightarrow\quad
 f\in p\,\mathcal D_d(M).
\]
This supplies a smooth definable quotient from formal divisibility;
definability of an already existing smooth quotient also follows
directly by taking the closure of its graph off the zero set of
the denominator.

For the constructed weights, formal-to-genuine solvability within
the original Denjoy--Carleman class fails even for the scalar equation
\[
 P(x,y)z-f(x,y)=0,
\]
linear in the unknown $z$, with polynomial coefficient $P$ vanishing
at the origin. It has a formal solution $\widehat g$ and a smooth
solution $g$, but no solution in $A_{M,2}$. Nevertheless,
\[
 g\in\mathcal D_2(M)\setminus A_{M,2};
\]
see Proposition~\ref{prop:definable-quotient}. Thus these examples
also exhibit strict inclusion of the fixed class in its smooth
definable enlargement.

Geometric finiteness results concern a different aspect of the
theory. Topological Noetherianity of quasianalytic zero sets, proved
by Bierstone--Milman \cite{BierstoneMilman2004}, and Nullstellensatz
statements such as \cite{AcquistapaceBrogliaNicoara2016} do not by
themselves imply algebraic Noetherianity of the local ring.
Related nonlinear regularity questions, including whether
membership of two consecutive powers of a function in an
Denjoy--Carleman class forces membership of the function itself,
are studied in \cite{NenningRainerSchindl2021}.

\subsection*{Our use of the appendices.} 
We have organized the paper so that the main text emphasizes the principal ideas and new mechanisms of the proofs. To this end, whenever possible, we have moved proofs based on adaptations of known methods to the appendices, retaining only some shorter proofs in the main text. We hope that this organization will help the reader identify the main novelties of our approach and distinguish them from the supporting, and 
well-known to specialists, arguments.

\subsection*{Acknowledgement} We are deeply grateful to Edward Bierstone, who introduced us to this problem and with whom we pursued earlier approaches to its solution. Our discussions with him, together with our earlier joint work on quasianalytic equations, composition, and loss of regularity under geometric transformations \cite{BelottoBiborskiBierstone2017,BelottoBierstoneChow2018,BelottoBierstoneKiro2022}, played an important role in shaping both the questions addressed here and our approach to them. Our understanding of these problems owes much to his mathematical insight and generosity.

\subsection*{AI assistance}
The research leading to this paper began in 2021, following the completion of \cite{BelottoBierstoneKiro2022}, and builds on ideas and constructions developed by the authors over the subsequent years. More recently, interactions with ChatGPT (OpenAI, Astra 6 Pro) accelerated this research. In addition to editorial assistance, ChatGPT contributed to the development of mathematical ideas and arguments, as detailed below.

The quantitative division theorem \ref{thm:positive} in Section \ref{sec:positive} grew out of an earlier proof by the authors for the logarithmic-class case of Theorem \ref{thm:intro-noetherianity}. Starting from this proof, ChatGPT analyzed the requirements for extending the argument, leading, through dialogue with the authors, to the identification of the condition $(G_\beta)$ and the development of a proof of the resulting generalization.

The negative results build on the authors' earlier strategy of using staircase weights with piecewise Gevrey growth. A counterexample in a sectorial quasianalytic Gevrey class, first obtained by ChatGPT and presented in Appendix~D, helped clarify how this approach could be implemented. The authors then introduced a polynomial approximation method for cutoff functions and, through further interactions with ChatGPT, developed a counterexample in a quasianalytic Denjoy--Carleman class. ChatGPT subsequently contributed to sharpening the construction and adapting it to the condition $(G_\beta)$, leading to the results presented in Sections~4--5.

Throughout these developments, the authors provided mathematical input and critically examined and refined the proposed arguments. All mathematical arguments in the final manuscript have been checked by the authors, who take full responsibility for the statements, proofs, and references.

\section{Weights and the complex divisor}\label{sec:weights}

We retain the standing assumptions of Definition~\ref{def:QDC} and
\eqref{eq:proper-weight}. The standard stability properties of the class
are recalled in Remark~\ref{rem:QDCbasic-properties}; see
\cite[\S\S2--4]{BierstoneMilman2004} and \cite{Thilliez2008}.
Here we record the weight estimates and branch geometry used in
Section~\ref{sec:positive}.

\subsection{Elementary weight calculus}
For a normalized logarithmically convex weight, set
$\mu_k=M_k/M_{k-1}$, $k\geq1$. Under \eqref{eq:proper-weight},
the sequence $\mu_k$ is nondecreasing and tends to infinity.
In particular, $h_M(t)>0$ and the infimum defining it is attained
for every $t>0$. The functions $h_M$ and $\omega_M$ in
\eqref{eq:associated} are the usual associated functions; see
\cite{Mandelbrojt1952} and \cite[\S7.2]{Rainer2022Survey}.
They are formed from the factorial-free weight $M$, which is
denoted by $m$ in the latter reference.

The function $\omega_M$ is defined for every $T>0$, whereas the
auxiliary notation $\Lambda_M=\log\omega_M$ is used only for
sufficiently large $T$. In terms of the standard associated function,
condition $(G_\beta)$ reads
\[
 \frac{\omega_M(T)}{\omega_M(T/R)}
 \geq \exp(cT^\beta)\qquad(T\gg1),
\]
for some fixed $R>1$ and $c>0$. Thus it measures growth between two
scales. The case $\beta=0$ is the usual associated-function
characterization of moderate growth; compare
\cite[Lemma~7.2 and footnote~8]{Rainer2022Survey}.
It permits absorption of fixed losses by changing scale.
For $\beta>0$, the stronger growth requirement will compensate
for the scale-dependent loss of continuation in
Section~\ref{sec:positive}. Neither of the following two lemmas
uses quasianalyticity.

\begin{lemma}\label{lem:weight-calculus}
The following assertions hold.
\begin{enumerate}
\item $M_jM_k\leq M_{j+k}$, and
\begin{equation}\label{eq:weight-duality}
 M_k=\sup_{T>0}T^k e^{-\omega_M(T)}.
\end{equation}
Moreover $h_M(t)=o(t^N)$ as $t\downarrow0$, for every fixed $N$.
\item Condition $(G_0)$ is equivalent to moderate growth:
\begin{equation}\label{eq:moderate-growth}
 M_{j+k}\leq H^{j+k}M_jM_k\qquad(j,k\geq0)
\end{equation}
for some $H\geq1$. Each $(G_\beta)$ implies $(G_{\beta'})$ when
$0\leq\beta'\leq\beta$.
\item Each $(G_\beta)$ is invariant under exponential equivalence
of weights: if
\[
 C^{-k}M_k\leq N_k\leq C^kM_k\qquad(k\geq0)
\]
for some $C\geq1$, then $M$ satisfies $(G_\beta)$ if and only if
$N$ does.
\end{enumerate}
\end{lemma}
\begin{proof}
The inequality in (i) follows by multiplying
$\mu_i\leq\mu_{j+i}$ for $1\leq i\leq k$.
If $\mu_k\leq T\leq\mu_{k+1}$, the term $T^k/M_k$ is maximal
among $T^j/M_j$, so equality is attained in
\eqref{eq:weight-duality}. For $k=0$, take $T\leq\mu_1$.
Finally, $h_M(t)\leq t^{N+1}M_{N+1}=o(t^N)$.

The monotonicity in $\beta$ is immediate from \eqref{eq:Gbeta}.
Under $(G_0)$, iteration gives $R_1>1$ and $C_0\geq0$ such that
\[
 \omega_M(T)\geq2\omega_M(T/R_1)-C_0\qquad(T>0).
\]
Extend the inequality from large $T$ by increasing $C_0$.
By \eqref{eq:weight-duality},
\[
 M_{j+k}\leq e^{C_0}\sup_{T>0}
 \bigl(T^je^{-\omega_M(T/R_1)}\bigr)
 \bigl(T^ke^{-\omega_M(T/R_1)}\bigr)
 \leq e^{C_0}R_1^{j+k}M_jM_k.
\]
Absorb the fixed factor when $j+k\geq1$. Conversely, taking
the infimum in
\[
 t^{j+k}M_{j+k}\leq(Ht)^jM_j(Ht)^kM_k
\]
gives $h_M(t)\leq h_M(Ht)^2$. Thus
$h_M(t)^{1/2}\leq h_M(Ht)$, which is $(G_0)$
(enlarge $H$ if necessary).

For (iii), the weight comparison gives
\[
 \Lambda_M(T/C)\leq\Lambda_N(T)\leq\Lambda_M(CT)
\]
for large $T$. If $M$ satisfies $(G_\beta)$ with constants $c,R$,
then
\begin{align*}
 \Lambda_N(T)-\Lambda_N\!\left(\frac{T}{C^2R}\right)
 &\geq\Lambda_M(T/C)-\Lambda_M(T/(CR))\\
 &\geq cC^{-\beta}T^\beta.
\end{align*}
Interchanging $M$ and $N$ proves the converse.
\end{proof}

The next lemma combines the standard absorption estimates for
fixed powers and polynomial factors; compare
\cite[\S1.2, inequalities~(6)--(7)]{Thilliez2002Sharp}.
For this assertion, moderate growth suffices.

\begin{lemma}[Absorbing fixed losses]\label{lem:weight-loss}
Assume moderate growth. For every $a\geq0$, $\theta>0$ and $B>0$,
there are $C,B'>0$ such that
\begin{equation}\label{eq:fixed-loss}
 t^{-a}h_M(Bt)^\theta\leq C h_M(B't)\qquad(0<t\ll1).
\end{equation}
\end{lemma}
\begin{proof}
Iteration of $h_M(t)\leq h_M(Ht)^2$ absorbs any fixed positive
power: choose $\ell$ such that $2^{-\ell}\leq\theta$ and use
$0<h_M\leq1$. For an integer $k\geq0$, moderate growth also gives
\[
 t^{-k}h_M(t)\leq\inf_{n\geq0}t^nM_{n+k}
 \leq H^k M_k h_M(Ht).
\]
Take an integer $k\geq a$, first absorb the power and then this
polynomial loss, and rescale $t$.
\end{proof}

\begin{remark}
We also record two elementary estimates used below.
For each fixed integer $s\geq0$, iteration of
\eqref{eq:standing-diff} gives
\begin{equation}\label{eq:fixed-shift}
 M_{k+s}\leq D^{sk+s(s+1)/2}M_k.
\end{equation}
If $N$ minimizes $t^kM_k$, then
\begin{equation}\label{eq:below-minimum}
 t^kM_k\leq1\qquad(0\leq k\leq N),
\end{equation}
since $k\mapsto\log(t^kM_k)$ is convex, starts at zero, and is
nonpositive at its minimum.
\end{remark}

\subsection{The branchwise separation exponent}\label{sec:geometry}
We use the branchwise separation exponent denoted by $d(\varphi)$ in
\cite[\S\S1.3--1.4]{Thilliez2002Sharp}.
Let $0\ne \varphi\in\R\{x,y\}$ vanish at the origin, and let
$\varphi_{\mathrm{red}}$ be the product of its distinct irreducible
analytic factors. Choose real linear coordinates in which the
lowest homogeneous part of $\varphi_{\mathrm{red}}$ does not vanish
on the $y$-axis. Weierstrass preparation writes
$\varphi_{\mathrm{red}}=UW$, with $U$ an analytic unit and
\begin{equation}\label{eq:adapted-weierstrass}
 W(x,y)=y^m+a_1(x)y^{m-1}+\cdots+a_m(x),\qquad
 a_j\in\R\{x\},\quad\ord_0a_j\geq j.
\end{equation}
Here $W$ is reduced. The coordinate choice is fixed throughout
the corresponding division argument.

By Newton--Puiseux, the roots of $W(z,\cdot)$ can be labelled by
holomorphic Puiseux branches $\lambda_\nu(z)$ on fixed sufficiently
small sectors about each real half-axis; see
\cite[\S4]{Teissier2007}. In this subsection, the index $\nu$
records both a root and the choice of half-axis, and the sign
$\pm$ below is the corresponding fixed sign.
On its chosen half-axis, a branch is real valued, or there is
a rational $q_\nu\geq1$ such that
\begin{equation}\label{eq:branch-imaginary}
 |\operatorname{Im}\lambda_\nu(\pm t)|
 \asymp t^{q_\nu}\qquad(t\downarrow0).
\end{equation}
We define
\begin{equation}\label{eq:branch-q}
 q(\varphi)=\max\bigl(\{1\}\cup\{q_\nu:\lambda_\nu
 \text{ is nonreal on its chosen half-axis}\}\bigr),
 \qquad \beta(\varphi)=1-\frac1{q(\varphi)}.
\end{equation}
For an analytic unit, set $q(\varphi)=1$.
The exponents in \eqref{eq:branch-imaginary} are the
regular-separation exponents of the corresponding nonreal branch
pieces from $\R^2$. Thus $q(\varphi)=d(\varphi)$ in Thilliez's notation.
For general background on separation exponents, see
\cite{BochnakRisler1975}.

The significance of $\beta(\varphi)$ is already visible in
\eqref{eq:branch-imaginary}. If a nonreal branch has height
$O(\sigma)$ at $\pm t$, then
\[
 t\lesssim\sigma^{1/q_\nu}\leq\sigma^{1/q(\varphi)}.
\]
Continuation along that branch through discs of radius comparable
to $\sigma$ can therefore require
$O(1+\sigma^{-\beta(\varphi)})$ steps.
This explains the exponent in the continuation argument of
Section~\ref{sec:positive}; the precise absorption of the loss
uses two scales, as in Lemma~\ref{pos:absorption}.

\begin{lemma}[Branch geometry]\label{lem:branch-geometry}
On the sectors above, the Puiseux branches satisfy
\begin{equation}\label{eq:branch-lipschitz}
 |\lambda_\nu(z)|\leq C|z|,\qquad
 |\lambda'_\nu(z)|\leq C.
\end{equation}
The first estimate holds for every root at every sufficiently
small complex $z$. Every nonreal branch satisfies
\eqref{eq:branch-imaginary} with $q_\nu\geq1$, and
\begin{equation}\label{eq:branch-monotone}
 |\operatorname{Im}\lambda_\nu(\pm t)|
 \leq C|\operatorname{Im}\lambda_\nu(\pm u)|
 \qquad(0<t\leq u\ll1).
\end{equation}
For some $c,A_0>0$, distinct roots on the same sector satisfy
\begin{equation}\label{eq:root-separation}
 |\lambda_i(z)-\lambda_j(z)|\geq c|z|^{A_0}.
\end{equation}
\end{lemma}
\begin{proof}
For $z\ne0$, the rescaled root $v=\lambda/z$ satisfies the monic
equation
\[
 v^m+\sum_{j=1}^m\frac{a_j(z)}{z^j}v^{m-j}=0.
\]
Its coefficients are bounded by \eqref{eq:adapted-weierstrass},
so its roots are uniformly bounded. The resulting estimate
$|\lambda_\nu(z)|\leq C|z|$ forces every nonzero leading exponent
in the convergent Puiseux expansions to be at least one.
Termwise differentiation gives the second estimate in
\eqref{eq:branch-lipschitz}. The imaginary part of each expansion
on its chosen half-axis is either identically zero or has a first
nonzero term. This proves \eqref{eq:branch-imaginary} and
\eqref{eq:branch-monotone}.

Since $W$ is reduced, the difference of two distinct Puiseux roots
is a nonzero convergent Puiseux series. Its leading term gives a
lower bound by a positive constant times a power of $|z|$.
There are finitely many pairs, so \eqref{eq:root-separation}
follows after shrinking the sectors and choosing common constants.
For $m=1$ the assertion is empty.
\end{proof}

Multiplicities and analytic units do not enter the definition.
Real-valued Puiseux branches impose no additional threshold.
This does not mean that a factor with real zeros may be discarded:
its roots can be real over one half-axis and nonreal over the
other. Nor is $q(\varphi)$ the maximum mutual contact order of the
complex branches. The following examples distinguish these points:
\begin{center}
\begin{tabular}{@{}lc@{}}
\toprule
Divisor $\varphi$ & $q(\varphi)$\\
\midrule
$y^2-x^{2k}$, $k\geq1$ & $1$\\
$y^2+x^{2k}$, $k\geq1$ & $k$\\
$y^2-x^3$ & $3/2$\\
$y^{2m}+x^{2n}$, $1\leq m<n$ & $n/m$\\
\bottomrule
\end{tabular}
\end{center}
For example, the two real branches of $y^2-x^{2k}$ have mutual
contact order $k$, whereas $q(\varphi)=1$. For $y^2-x^3$, the negative
half-axis gives the nonreal roots and hence the exponent $3/2$.

The exponent for the whole complex zero set can also differ from
\eqref{eq:branch-q}. In Thilliez's example
\[
 P=y\bigl((y-x^2)^2+x^8\bigr),
\]
the nonreal roots are $x^2\pm ix^4$, so $q(P)=4$, while the
regular-separation exponent of the entire complex zero set
$Z_P\subset\C^2$ from the real plane is $2$
\cite[Example~1.5(iv)]{Thilliez2002Sharp}.
Along the nonreal branch $z_t=(t,t^2+it^4)$,
\[
 \dist(z_t,\R^2)\asymp t^4,\qquad
 \dist(z_t,Z_P\cap\R^2)\asymp t^2.
\]
The whole-set separation inequality measures distance to the
real zero set $Z_P\cap\R^2=\{y=0\}$, whereas each of these
nonreal branches meets the real plane only at the origin.

\subsection{Logarithmic examples}
\begin{proposition}\label{prop:log-weights}
Let $M_k=(\log(k+e))^{ak}$, where $0<a\leq1$.
Then the standing assumptions hold, the class is proper, and
$(G_\beta)$ holds for every $0\leq\beta<1$.
\end{proposition}
\begin{proof}
The function $x\mapsto x\log\log(x+e)$ is increasing on
$[0,\infty)$. If $u=\log(x+e)$, its second derivative is
\[
 \frac{x(u-1)+2eu}{(x+e)^2u^2}>0.
\]
This proves monotonicity and logarithmic convexity.
Also $M_k^{1/k}\to\infty$, and the mean value theorem gives
$M_{k+1}/M_k\sim(\log k)^a$. This implies differentiation
stability and the divergence in \eqref{eq:standing-qa}
for $a\leq1$.

Put $v=T^{1/a}$. Then
\[
 \omega_M(T)=\sup_{k\geq0}ak\log\frac{v}{\log(k+e)}.
\]
A positive term has $k<e^v$ and is at most $ae^v\log v$.
Taking $k=\lfloor e^{v-1}\rfloor$ gives a lower bound
$c_a e^v/v$ for large $v$. Therefore
\[
 \Lambda_M(T)=T^{1/a}+O(\log T).
\]
For each fixed $R>1$, subtraction gives
\[
 \Lambda_M(T)-\Lambda_M(T/R)
 =(1-R^{-1/a})T^{1/a}+O(\log T),
\]
which implies every $(G_\beta)$ in the stated range.
\end{proof}

By Theorem~\ref{thm:plane}, these logarithmic weights therefore
give regular Noetherian plane local rings $A_{M,2}$ of dimension
two. In contrast, the staircase weights of
Section~\ref{sec:staircase} satisfy $(G_\beta)$ precisely for
$0\leq\beta<\beta_0$, where $\beta_0\in(0,1)$ is prescribed.
\section{Division by analytic plane germs}\label{sec:positive}

We prove Theorem~\ref{thm:positive} by constructing holomorphic
approximants to the prescribed smooth quotient. We first recall the
approximation statements that will be used at the beginning and at the
end of the proof. Throughout this section the standing assumptions on
$M$ hold, $M_k^{1/k}\to\infty$, and $M$ has moderate growth. The latter
assumption follows from every condition $(G_\beta)$ considered in the
theorem, by Lemma~\ref{lem:weight-calculus}.

\subsection{Holomorphic approximation}

Beurling described quasianalytic classes through the balance between
the accuracy of approximation and the size of holomorphic approximants,
on both fixed domains and domains shrinking towards a real arc
\cite[Chapter~III, \S1, classes~IV and~V, pp.~397--400]{Beurling1989}.
Dynkin's almost analytic description of Denjoy--Carleman classes gives
a way to construct such approximants: extend the function to a complex
neighbourhood with a small $\bar\partial$-derivative, then correct that
derivative by a Cauchy--Green operator
\cite{Dynkin1980,Dynkin1993}. Conversely, small differences between
successive holomorphic approximants give derivative estimates by
Cauchy's inequalities. Holomorphic approximation characterizations on
intervals are given in
\cite[Proposition~3.3.2 and Remark~3.3.3]{Thilliez2020Powers} and
\cite[Theorem~5.3]{NenningRainerSchindl2021}.

We record the two forms needed here. The first provides bounded
approximants on a tube of a fixed rectangle; the second allows local
approximants whose norms grow as a fixed power of the inverse radius.
Their proofs, including the several-variable extension input, are in
Appendix~\ref{app:holomorphic}.

For $d=1,2$ write
\[
 D_{\C}^{d}(p,r)=\{z\in\C^d:\max_j|z_j-p_j|<r\},\qquad
 D_{\R}^{d}(p,r)=D_{\C}^{d}(p,r)\cap\R^d,
\]
and, for a compact real rectangle $K\subset\R^2$, put
\[
 T_r(K)=\bigcup_{p\in K}D_{\C}^{2}(p,r).
\]

\begin{proposition}[Bounded holomorphic approximants]\label{prop:hol-approx}
Let $u\in\CM(U)$, where $U\subset\R^2$ is open, and let $K$ be a
compact real rectangle contained in a ball $B\Subset U$. For every
fixed $T>0$ there are $A,C,s_0>0$ and holomorphic functions $u_s$ on
$T_{Ts}(K)$, $0<s<s_0$, such that
\begin{equation}\label{eq:bounded-hol-approx}
 \sup_{T_{Ts}(K)}|u_s|\le C,\qquad
 \sup_K|u_s-u|\le C h_M(As).
\end{equation}
The constants may depend on $T$, but not on $s$.
\end{proposition}

This is a shrinking-tube consequence of Dynkin's almost analytic
extension theorem, in the several-variable formulation of
\cite[Theorem~1.2]{FurdosNenningRainerSchindl2020}.
The construction uses the same extension and Cauchy--Green correction
as the one-variable forward implication in
\cite[Proposition~3.3.2 and Remark~3.3.3]{Thilliez2020Powers}.
We prove the stated two-variable form in
Appendix~\ref{app:holomorphic}.

\begin{proposition}[A local approximation criterion]
\label{prop:approx-characterization}
Let $u$ be smooth on a neighbourhood of a compact set $K\subset\R^2$.
Suppose that there are $\kappa,B,C,\sigma_0>0$ and $a\ge0$ such that,
for every $p\in K$ and $0<\sigma<\sigma_0$, there is a holomorphic
function $U_{p,\sigma}$ on $D_{\C}^{2}(p,\kappa\sigma)$ satisfying
\begin{equation}\label{eq:local-approx-hypothesis}
 \sup_{D_{\C}^{2}(p,\kappa\sigma)}|U_{p,\sigma}|\le C\sigma^{-a},
 \qquad
 \sup_{D_{\R}^{2}(p,\kappa\sigma)}|U_{p,\sigma}-u|
 \le C h_M(B\sigma).
\end{equation}
Assume that these real squares lie in the domain of $u$. Then there
are $A_1,B_1>0$ such that
\[
 \sup_{p\in K}|\partial^\alpha u(p)|
 \le A_1B_1^{|\alpha|}\alpha!M_{|\alpha|}
 \qquad(\alpha\in\N^2).
\]
\end{proposition}

This is a local two-variable variant of the converse approximation
theorems in \cite[Proposition~3.3.2]{Thilliez2020Powers} and
\cite[Theorem~5.3(ii)]{NenningRainerSchindl2021}.
Here the approximants may depend on the centre, and their complex
norms may grow polynomially as the radius shrinks. The proof of
this form is given in Appendix~\ref{app:holomorphic}.

The constants in the converse must be uniform in $p$ and $\sigma$.
The functions $U_{p,\sigma}$ need not agree on overlaps. In particular,
uniform local approximants on a fixed real rectangle prove membership
in $\CM$ on its interior. Moderate growth is used in the converse to
absorb the fixed powers arising from real-to-complex smallness.

\subsection{The division strategy}

We first treat a reduced divisor. In the adapted coordinates of
Section~\ref{sec:geometry}, let $W$ be a reduced Weierstrass polynomial
of degree $m\ge1$, and suppose that
\begin{equation}\label{pos:identity}
 f=Wg,\qquad f\in\CM,\quad g\in C^\infty
\end{equation}
on a fixed neighbourhood of the origin. Put $q=q(W)$ and
$\beta_W=1-1/q$. The hypothesis $q\le(1-\beta)^{-1}$ gives
$\beta_W\le\beta$; hence it is enough to assume $(G_{\beta_W})$.
Write $D_r(p)=D_{\C}^{2}(p,r)$ and
$D_r^{\R}(p)=D_{\R}^{2}(p,r)$, and omit $p$ when $p=0$.

Choose a fixed real rectangle $K$ around the origin, contained in a
ball on which \eqref{pos:identity} holds. Proposition~\ref{prop:hol-approx}
provides holomorphic functions $F_s$ on $T_{Ts}(K)$ such that
\begin{equation}\label{pos:approximants}
 \|F_s\|_{T_{Ts}(K)}\le C,\qquad
 \|F_s-f\|_K\le C h_M(As).
\end{equation}
The tube factor $T$ will be fixed from the geometry before choosing
the relation between the two small scales.

Our main task is the following local division assertion.
\begin{proposition}[Local approximate division]\label{pos:local-division}
There are $r,\kappa,A,C,\sigma_0>0$ such that, for every
$p\in[-r,r]^2$ and $0<\sigma<\sigma_0$, there is a holomorphic
function $G_{p,\sigma}$ on $D_{2\kappa\sigma}(p)$ satisfying
\begin{equation}\label{pos:local-division-estimates}
 \|G_{p,\sigma}\|_{D_{2\kappa\sigma}(p)}\le C\sigma^{-m},\qquad
 \|f-WG_{p,\sigma}\|_{D_{2\kappa\sigma}^{\R}(p)}
 \le C h_M(A\sigma).
\end{equation}
\end{proposition}

The integral construction in the proof is the classical Cauchy
construction for holomorphic Weierstrass division; see
\cite[Theorem~6.2.5]{LeblSCV}. Its quantitative consequences are
recorded in Lemma~\ref{app:cauchy-division}. Below we choose the
contours to respect both the approximation tube and the region
where smallness on the divisor is available.

Here $\sigma$ is the radius needed for the local approximants to $g$.
We approximate $f$ at the larger scale $s=L\sigma$, where $L\ge1$
will be fixed independently of $\sigma$. Formal divisibility at the
origin first makes $F_s$ small on the central part of $\{W=0\}$.
Propagation along its Puiseux branches then gives smallness on the
parts of the divisor within distance $O(\sigma)$ of the real plane.
Condition $(G_{\beta_W})$ compensates for the continuation loss.
A suitable Cauchy integral produces $G_{p,\sigma}$ with
\eqref{pos:local-division-estimates}. Finally, the identity
\[
 W(g-G_{p,\sigma})=f-WG_{p,\sigma}
\]
allows us to apply a finite-difference estimate and remove the factor
$W$. The resulting
approximants to $g$ meet Proposition~\ref{prop:approx-characterization}.

The branchwise geometry is the same geometry used in
\cite[\S\S1.3--1.4]{Thilliez2002Sharp}. The additional estimate needed
here controls the loss from continuation at the two scales $s$ and
$\sigma$ while remaining in the original class.

\subsection{Formal divisibility and central smallness}

The identity \eqref{pos:identity} implies
$\widehat f=W\widehat g$. The formal division estimate
\cite[Theorem~1]{Thilliez2002Formal} gives
\[
 \widehat g=\sum_{\alpha\in\N^2}g_\alpha X^\alpha,
 \qquad |g_\alpha|\le A_0B^{|\alpha|}M_{|\alpha|}.
\]
An elementary proof in the present coordinates is included in
Lemma~\ref{pos:formal-coefficients}. We use this estimate only at the
origin; it does not yet give derivative bounds for $g$ nearby.

\begin{lemma}[Central smallness]\label{pos:central}
Fix $K_0>1$, and choose the tube in \eqref{pos:approximants} to contain
$D_{2K_0s}$. There are polynomials $Q_s$ and constants $A_{\mathrm{cen}},C>0$ such that
\begin{equation}\label{pos:central-bounds}
 \|Q_s\|_{D_{2K_0s}}\le C,\qquad
 \|F_s-WQ_s\|_{D_{K_0s}}\le C h_M(A_{\mathrm{cen}}s).
\end{equation}
Consequently, $|F_s|\le C h_M(A_{\mathrm{cen}}s)$ on
$\{W=0\}\cap D_{K_0s}$.
\end{lemma}

The real Taylor-remainder argument is the analogue, for analytic
$W$, of \cite[Lemma~4, inequalities~(20)--(21)]{Thilliez2002Formal}.
We optimize the truncation order and then pass from real to complex
smallness; the full proof is in
Appendix~\ref{app:central-proof}.

The construction is $Q_s=T_{N(s)-1}\widehat g$, where $N(s)$
minimizes $(C_*s)^NM_N$ for a sufficiently large fixed $C_*$.
The formal coefficient estimate makes $Q_s$ uniformly bounded
on a polydisc of radius comparable to $s$. Formal divisibility gives
$T_{N-1}f=T_{N-1}(WQ_s)$, so Taylor's formula and the analytic
tail estimate make $f-WQ_s$ small on the real square. Adding
$F_s-f$ and applying the real-to-complex estimate gives
\eqref{pos:central-bounds}.

Thus it is $F_s-WQ_s$ that is small throughout the central polydisc.
Restricting to $W=0$ gives the initial smallness of $F_s$ that will
be propagated along the branches.

\subsection{Propagation along the complex divisor}

We use the Puiseux branch geometry recalled in
\cite[\S\S1.3--1.4]{Thilliez2002Sharp}, in the form proved in
Lemma~\ref{lem:branch-geometry}. On sectors about the real half-axes,
that lemma
gives
\begin{equation}\label{pos:geometry-used}
 |\lambda_\nu(z)|\le C|z|,\qquad |\lambda'_\nu(z)|\le C.
\end{equation}
On either half-axis, a branch is real valued or satisfies
$|\operatorname{Im}\lambda_\nu(\pm t)|\asymp t^{q_\nu}$,
where $1\le q_\nu\le q$. In the latter case,
\begin{equation}\label{pos:height-monotone}
 |\operatorname{Im}\lambda_\nu(\pm t)|
 \le C|\operatorname{Im}\lambda_\nu(\pm u)|
 \qquad(0<t\le u\ll1).
\end{equation}

For a nonreal branch with height at most $4\sigma$ at $a\in\R$,
we have $|a|\le C\sigma^{1/q}$. Propagation from a point at distance
$s=L\sigma$ from the origin to $a$, through discs of radius
comparable to $s$, therefore requires at most
$C(1+L^{-1}\sigma^{-\beta_W})$ steps. Each step replaces the
smallness by a fixed fractional power. The next lemma states how
$(G_{\beta_W})$ compensates for the resulting loss of smallness.

\begin{lemma}[Absorbing the continuation loss]\label{pos:absorption}
Assume $(G_\gamma)$, where $0\le\gamma<1$. For every $a,D_1>0$
there are $L\ge1$ and $b>0$ such that
\begin{equation}\label{pos:absorption-estimate}
 e^{-D_1(1+L^{-1}\sigma^{-\gamma})}
 \omega_M\!\left(\frac{a}{L\sigma}\right)
 \ge\omega_M\!\left(\frac{b}{\sigma}\right)
 \qquad(0<\sigma\ll1).
\end{equation}
\end{lemma}

\begin{proof}
For $\gamma>0$, write $(G_\gamma)$ as
$\Lambda_M(T)-\Lambda_M(T/R)\ge cT^\gamma$.
Choose $L$ so that $ca^\gamma L^{1-\gamma}>2D_1$, and put
$b=a/(RL)$. For sufficiently small $\sigma$,
\[
 \Lambda_M\!\left(\frac{a}{L\sigma}\right)
 -\Lambda_M\!\left(\frac{b}{\sigma}\right)
 \ge ca^\gamma L^{-\gamma}\sigma^{-\gamma}
 \ge D_1(1+L^{-1}\sigma^{-\gamma}).
\]
Exponentiation gives \eqref{pos:absorption-estimate}. For $\gamma=0$,
take $L=1$ and iterate $(G_0)$ a fixed number $m$ of times with
$mc\ge2D_1$; then $b=a/R^m$ works.
\end{proof}

\begin{lemma}[Smallness on branches near the real plane]
\label{pos:branch-smallness}
There are fixed $L\ge1$ and $r,\eta,A,C,\sigma_0>0$ such that,
with $s=L\sigma$, the conditions
\[
 a\in\R,\qquad 0<\sigma<\sigma_0,\qquad
 \sigma\le|a|\le r,\qquad
 |\operatorname{Im}\lambda_\nu(a)|\le4\sigma
\]
imply that $\lambda_\nu$ is holomorphic on $|z-a|<\eta\sigma$ and
\begin{equation}\label{pos:branch-smallness-estimate}
 \sup_{|z-a|<\eta\sigma}|F_s(z,\lambda_\nu(z))|
 \le C h_M(A\sigma).
\end{equation}
\end{lemma}

The analytic propagation tool is the standard iteration of
Hadamard's three-circles theorem, recorded in
Lemma~\ref{lem:disc-chain}; compare
\cite[Theorem~3.1 and \S7]{Trefethen2020Continuation}.
The branch geometry determines the number of iterations, and
Lemma~\ref{pos:absorption} then compensates for their loss.

\begin{proof}
Choose the geometric constants $\eta,K_0,T$ before $L$, as
verified in Appendix~\ref{app:branch-geometric-details}. Thus the
graph discs used below lie in the fixed approximation tube,
and their uniform bounds have constants independent of $L$.
Put $H(z)=F_s(z,\lambda_\nu(z))$.
For $\sigma\le|a|\le s$, central smallness applies directly.
For a real branch with $|a|\ge s$, the identity
$f(t,\lambda_\nu(t))=0$ gives real smallness of $H$;
the real-to-complex estimate then
gives complex smallness. The details of these two cases are
also in Appendix~\ref{app:branch-geometric-details}.

For a nonreal branch with $|a|\ge s$, its height gives
\[
 |a|\le C\sigma^{1/q_\nu}\le C\sigma^{1/q}.
\]
Join $\operatorname{sgn}(a)s$ to $a$ by real centres separated by
at most $\eta s$. Along this interval \eqref{pos:height-monotone}
bounds the central heights by $C\sigma$. Thus $H$ is holomorphic
and uniformly bounded on the chain of discs of radius $4\eta s$.
The number $N$ of steps satisfies
\begin{equation}\label{pos:chain-length}
 N\le C\left(1+\frac{|a|}{s}\right)
 \le C(1+L^{-1}\sigma^{-\beta_W}).
\end{equation}
The first disc of radius $\eta s$ has its graph in $D_{K_0s}$,
so Lemma~\ref{pos:central} bounds $H$ there by $C h_M(A_{\mathrm{cen}}s)$.
Lemma~\ref{lem:disc-chain} gives on the last disc a bound of the form
\begin{equation}\label{pos:chain-result}
 C\exp\left\{-e^{-D_2(1+L^{-1}\sigma^{-\beta_W})}
       \omega_M\!\left(\frac{a_0}{L\sigma}\right)\right\}.
\end{equation}
Here $a_0=A_{\mathrm{cen}}^{-1}$ and $D_2>0$ are independent of $L$: the
geometric choices and $T$ were fixed first. Constant factors in
the initial smallness contribute only a fixed prefactor after
raising to an exponent in $(0,1]$.

Apply Lemma~\ref{pos:absorption} with $\gamma=\beta_W$ to choose
$L$ once and for all. The last bound is then at most
$C h_M(\sigma/b)$. Its disc has radius $\eta s\ge\eta\sigma$.
Enlarging $A$ also covers the central and real-branch cases.
\end{proof}

\begin{remark}\label{pos:endpoint-remark}
The continuation loss in \eqref{pos:chain-length} has coefficient
$L^{-1}$, whereas the gain from $(G_{\beta_W})$ has coefficient
$L^{-\beta_W}$. Since $\beta_W<1$, their ratio is improved by
increasing $L$, which remains fixed as $\sigma\downarrow0$.
This proves the endpoint $q(W)=(1-\beta)^{-1}$ as well as the
strict inequality case.
\end{remark}

\subsection{Construction of approximate quotients}

The natural candidate is $F_s/W$. It need not be holomorphic near
$p$: although $f=Wg$ on the real plane, $F_s$ need not vanish on
the complex zeros of $W$. The preceding estimates say that its
values there are small. We therefore use a Cauchy integral of
$F_s/W$ to construct a holomorphic function inside a contour.
The residues of the poles give an explicitly controlled remainder.

Lemma~\ref{app:cauchy-division} supplies the standard integral and
residue estimates. The choice of contour is the essential point. It must enclose a
neighbourhood of $p$ of radius comparable to $\sigma$, remain at
distance comparable to $\sigma$ from every root, and lie in the
tube where $F_s$ is defined. Moreover, every enclosed root must
belong to the part of the divisor on which we have smallness.
Near the origin roots may collide, so we use the central estimate
in place of separate residue bounds.

\begin{proof}[Proof of Proposition~\ref{pos:local-division}]
Use $F_s$ with $s=L\sigma$, where $L$ is fixed by
Lemma~\ref{pos:branch-smallness}. All constants multiplying $s$
may now be replaced by fixed constants multiplying $\sigma$.
We construct $G$ on a polydisc of radius $c\sigma$ about
$p=(a,b)$, uniformly in $p,\sigma$.

\smallskip\noindent\emph{Near the origin: $|a|\le\sigma$.}
The root bound of Lemma~\ref{lem:branch-geometry} places every
root of $W(z,\cdot)$ in $|w|\le C_0\sigma$ for $|z|<2\sigma$.
Choose a fixed $R_0>C_0+4$. If $|b|\le(R_0-2)\sigma$, use the
circle $|\zeta|=R_0\sigma$ and set
\begin{equation}\label{pos:central-integral}
 G(z,w)=\frac{1}{2\pi i}\int_{|\zeta|=R_0\sigma}
 \frac{F_s(z,\zeta)}{W(z,\zeta)(\zeta-w)}\,d\zeta,
 \qquad |z|<2\sigma,\quad |w|<(R_0-1)\sigma.
\end{equation}
The circle encloses every root and stays a fixed multiple of
$\sigma$ from them. Lemma~\ref{app:cauchy-division} makes $G$
holomorphic with $\|G\|\le C\sigma^{-m}$.

On this circle and its interior, Lemma~\ref{pos:central} gives
$F_s=WQ_s+H_s$, with $\|H_s\|\le C h_M(A\sigma)$.
The estimate in Lemma~\ref{app:cauchy-division} for this
decomposition, which also applies when roots collide, gives
\[
 |F_s-WG|\le C(1+\sigma^{-m})h_M(A\sigma).
\]
Absorb the fixed power by Lemma~\ref{lem:weight-loss} and add
$f-F_s$ on the real part. This proves the required estimates on
a uniform $c\sigma$-polydisc about $p$.

If $|b|>(R_0-2)\sigma$, a sufficiently small $c\sigma$-polydisc
about $p$ stays at distance at least $c_1\sigma$ from every root.
Here take $G=F_s/W$ directly. Then $\|G\|\le C\sigma^{-m}$,
and $f-WG=f-F_s$ on the real part.

\smallskip\noindent\emph{Away from the origin: $|a|\ge\sigma$.}
On $|z-a|<c\sigma$, with $c>0$ sufficiently small, the roots are
holomorphic and distinct. Lemma~\ref{lem:branch-geometry} gives
\begin{equation}\label{pos:root-separation}
 |\lambda_i(z)-\lambda_j(z)|\ge c_2\sigma^{A_0}\qquad(i\ne j),
\end{equation}
because $|z|\ge\sigma/2$.

Choose a rectangle $\Gamma$ centred at $b\in\R$ in the $w$-plane.
Choose its horizontal half-width $u\in[\sigma,2\sigma]$ at
distance at least $\sigma_m\sigma$ from every number
$|\operatorname{Re}\lambda_i(a)-b|$, and its vertical half-width
$v\in[\sigma,2\sigma]$ at the same distance from every number
$|\operatorname{Im}\lambda_i(a)|$.
These choices are possible with $\sigma_m=1/(4(m+1))$, since
the excluded lengths in either interval total at most
$2m\sigma_m\sigma<\sigma$. Every root at $a$ is therefore
separated from $\Gamma$ by at least $\sigma_m\sigma$.
The derivative bound in \eqref{pos:geometry-used} preserves
separation by $\sigma_m\sigma/2$ on $|z-a|<c\sigma$, after
decreasing $c$. In particular the set of enclosed roots is fixed.

For $|w-b|<\sigma/2$ define
\begin{equation}\label{pos:local-integral}
 G(z,w)=\frac{1}{2\pi i}\int_\Gamma
 \frac{F_s(z,\zeta)}{W(z,\zeta)(\zeta-w)}\,d\zeta.
\end{equation}
The contour lies in the chosen approximation tube. It has
length $O(\sigma)$ and stays at distance comparable to
$\sigma$ from the roots and the target $w$-disc.
Lemma~\ref{app:cauchy-division} again gives holomorphicity and
$\|G\|\le C\sigma^{-m}$.

The same lemma expresses $F_s-WG$ through the values of $F_s$
at the enclosed roots. Each such root has
$|\operatorname{Im}\lambda_i(a)|<2\sigma$, precisely because
of the chosen vertical half-width. Lemma~\ref{pos:branch-smallness}
therefore bounds these values by $C h_M(A\sigma)$ on the
parameter disc, after a further fixed decrease of $c$.
Together with \eqref{pos:root-separation}, the quantitative
remainder estimate in Lemma~\ref{app:cauchy-division} gives
\[
 |F_s-WG|\le C\sigma^{-A_0(m-1)}h_M(A\sigma).
\]
Absorb this fixed loss and add $f-F_s$ on the real square.
This proves \eqref{pos:local-division-estimates}; choose
$2\kappa\le c$ to obtain the stated polydisc.
\end{proof}

\subsection{Recovery of the smooth quotient}

The proposition makes $W(g-G_{p,\sigma})$ small. We cannot divide
this error pointwise by $W$ near its real zeros. The following
elementary estimate uses a first derivative bound to transfer
smallness from a nearby point where the monic polynomial is large.

\begin{lemma}[Removal of a small factor]\label{pos:small-product}
Let $W(x,y)$ be monic of degree $m\ge1$ in $y$. Suppose that on
a real square of radius $r$, a $C^1$ function $v$ satisfies
\[
 |Wv|\le\varepsilon,\qquad |\partial_yv|\le L_0,
 \qquad 0<\varepsilon^{1/(m+1)}\le r/2.
\]
On the concentric square of radius $r/2$ one has
\begin{equation}\label{pos:small-product-bound}
 |v|\le C_m(1+L_0)\varepsilon^{1/(m+1)}.
\end{equation}
The assertion also holds for complex-valued $v$ and coefficients.
\end{lemma}

This elementary estimate follows from the $m$th finite-difference
identity. Its proof is included in
Appendix~\ref{app:small-product-proof}.

\begin{proof}[Proof of Theorem~\ref{thm:positive}]
Multiplication by an analytic function, division by an analytic
unit, and real linear changes of coordinates preserve $\CM$.
In particular the theorem is immediate for an analytic unit.

For a reduced divisor vanishing at the origin, write $\varphi=VW$ by
analytic Weierstrass preparation in adapted coordinates, with $V$
an analytic unit. Apply Proposition~\ref{pos:local-division} to
$u=f/V=Wg$. Cauchy's estimate gives
\[
 |\partial_yG_{p,\sigma}|\le C\sigma^{-m-1}
 \quad\hbox{on }D_{\kappa\sigma}(p).
\]
Since $g$ is smooth on a fixed neighbourhood, the same bound,
after enlarging $C$, holds for $\partial_y(g-G_{p,\sigma})$ on
the real square. Apply Lemma~\ref{pos:small-product} to
\[
 v=g-G_{p,\sigma},\qquad Wv=u-WG_{p,\sigma},
 \qquad \varepsilon=C h_M(A\sigma),\quad r=\kappa\sigma.
\]
Its smallness condition holds uniformly for small $\sigma$,
since $h_M(A\sigma)$ is smaller than every fixed power of
$\sigma$. We obtain
\[
 \|g-G_{p,\sigma}\|_{D_{\kappa\sigma/2}^{\R}(p)}
 \le C\sigma^{-m-1}h_M(A\sigma)^{1/(m+1)}
 \le C' h_M(A'\sigma),
\]
where the last step is Lemma~\ref{lem:weight-loss}. Restricting
the holomorphic functions to $D_{\kappa\sigma/2}(p)$ gives the
uniform local approximants required by
Proposition~\ref{prop:approx-characterization}. Thus $g\in\CM$
near the origin. No gluing of the $G_{p,\sigma}$ is needed: all
of them approximate the same prescribed smooth $g$, and the
criterion recovers its derivatives at each centre.

Finally factor an arbitrary nonzero analytic divisor as
$\varphi=V\prod_i\varphi_i^{e_i}$, where $V$ is a unit and the $\varphi_i$ are
distinct irreducible real analytic factors. Divide successively
by these factors, with their multiplicities. At every step the
prescribed quotient is smooth, being the product of the remaining
factors with the original $g$. Coordinates adapted to the reduced
product are adapted to each factor, and the branch exponents of
each factor occur among those of that product. Thus
$\beta(\varphi_i)\le\beta(\varphi)\le\beta$, and $(G_\beta)$ supplies the
required hypothesis at every step.
\end{proof}
\section{A staircase weight with a prescribed threshold}
\label{sec:staircase}

Fix $0<\beta_0<1$ and put $\alpha=1/\beta_0>1$.  We construct one
proper quasianalytic weight of moderate growth for which $(G_\beta)$
holds precisely when $0\leq\beta<\beta_0$.  Recall our normalization:
$|\partial^\gamma f|\leq AB^{|\gamma|}|\gamma|!M_{|\gamma|}$, and
\[
 \omega_M(T)=\sup_{k\geq0}\log\frac{T^k}{M_k},\qquad
 \Lambda_M(T)=\log\omega_M(T),\qquad
 (G_\beta):\ \Lambda_M(T)-\Lambda_M(T/R)\geq cT^\beta
\]
for some $R>1$, $c>0$, and all sufficiently large $T$.
The positive division theorem then applies to contact exponents
$q<(1-\beta_0)^{-1}$; Section~\ref{sec:counterexample} constructs
failures for $P=y^{2m}+x^{2n}$ when $q=n/m\geq(1-\beta_0)^{-1}$.
The weight is fixed before the divisor is chosen.  These are
counterexamples for the constructed weights, not a converse for an
arbitrary weight.

\subsection{Construction of the staircase}
We alternate growth intervals of the quotient sequence
$\mu_k=M_k/M_{k-1}$ with constant intervals.  The growth profile
sets the threshold, long plateaux ensure quasianalyticity, and long
growth intervals accommodate the polynomial degrees needed later.
For $t\geq1$, set
\begin{equation}\label{stair:LQ}
 L(t)=\log(t+e^e),\qquad
 Q(t)=L(t)^\alpha\bigl(\log L(t)\bigr)^{2\alpha}.
\end{equation}
Write $E_\ell(t)=\exp^{\circ\ell}(t)$, with $E_0(t)=t$, and prescribe
\begin{equation}\label{stair:breakpoints}
 r_j=\left\lceil E_{6j}(2)\right\rceil,\qquad
 s_j=\left\lceil E_{6j+4}(2)\right\rceil\qquad(j\geq1).
\end{equation}
Thus $r_j<s_j<r_{j+1}$.  These are derivative orders: $(r_j,s_j]$
is a growth interval and $(s_j,r_{j+1}]$ is a plateau.  Put
\begin{equation}\label{stair:cj}
 c_j=\frac1{Q(r_j)}\prod_{i<j}\frac{Q(s_i)}{Q(r_i)}.
\end{equation}
The matching constants give the quotient sequence
\begin{equation}\label{stair:piecewise}
 \mu_k=\begin{cases}
  1,&k\leq r_1,\\
  c_jQ(k),&r_j<k\leq s_j,\\
  c_jQ(s_j),&s_j<k\leq r_{j+1},
 \end{cases}
\end{equation}
and hence the explicit weight
\begin{equation}\label{stair:weight}
 M_0=1,\qquad M_k=\prod_{i=1}^k\mu_i.
\end{equation}
In the coordinate $u=\log Q(k)$, the graph of $\log\mu_k$ has slope
one on growth intervals and zero on plateaux; see
Figure~\ref{fig:staircase}.

\begin{figure}[htbp]
\centering
\begin{tikzpicture}[x=0.96cm,y=0.67cm,>=Latex,font=\small]
 \fill[blue!5] (1,0) rectangle (2.8,5.9);
 \fill[blue!5] (4.6,0) rectangle (6.5,5.9);
 \fill[blue!5] (8.0,0) rectangle (9.8,5.9);
 \draw[->] (0,0)--(11.2,0) node[above left] {$u=\log Q(k)$};
 \draw[->] (0,0)--(0,6.2) node[above right] {$\log\mu_k$};
 \foreach \x/\lab in {1/{\log Q(r_1)},2.8/{\log Q(s_1)},4.6/{\log Q(r_2)},6.5/{\log Q(s_2)},8.0/{\log Q(r_3)},9.8/{\log Q(s_3)}}
 {\draw (\x,0.08)--(\x,-0.08);\node[below,font=\scriptsize] at (\x,-0.10) {$\lab$};}
 \draw[densely dashed,gray] (2.8,0)--(2.8,1.8);
 \draw[densely dashed,gray] (4.6,0)--(4.6,1.8);
 \draw[densely dashed,gray] (6.5,0)--(6.5,3.7);
 \draw[densely dashed,gray] (8.0,0)--(8.0,3.7);
 \draw[densely dashed,gray] (9.8,0)--(9.8,5.5);
 \draw[very thick,blue!55!black] (0,0)--(1,0)--(2.8,1.8)--(4.6,1.8)--(6.5,3.7)--(8,3.7)--(9.8,5.5)--(10.7,5.5);
 \node[above,font=\small] at (3.7,1.8) {plateau};
 \node[above,rotate=34,font=\small] at (5.7,3.15) {growth};
 \node[above,font=\small] at (7.25,3.7) {plateau};
\end{tikzpicture}
\caption{The graph of $\log\mu_k$ against $u=\log Q(k)$.
The interval lengths are schematic; the breakpoints are prescribed
by \eqref{stair:breakpoints}.}
\label{fig:staircase}
\end{figure}

\begin{samepage}
There are two spacing requirements.  The second also supplies the
endpoint estimate used to exclude $(G_{\beta_0})$.
\begin{lemma}\label{stair:separation}
As $j\to\infty$,
\begin{equation}\label{stair:plateau-room}
 \frac{\log((r_{j+1}+1)/(s_j+1))}{Q(s_j)}\longrightarrow\infty,
\end{equation}
and, for every fixed $C>0$,
\begin{equation}\label{stair:degree-room}
 \frac{E_3(Cr_j)}{s_j}\longrightarrow0.
\end{equation}
In particular,
\begin{equation}\label{stair:endpoint-room}
 \frac{Q(r_j)^{1/\alpha}}{(\log L(s_j))^2}\longrightarrow0.
\end{equation}
Moreover, $j\leq\log r_j$ for every $j\geq1$.
\end{lemma}
\end{samepage}
The verification is in Appendix~\ref{app:stair-details}.  The explicit
towers avoid any later choice of breakpoints depending on the divisor.

\subsection{Regularity and the exact growth threshold}
\begin{proposition}\label{stair:regularity}
The weight \eqref{stair:weight} is increasing and logarithmically
convex, is stable under differentiation, and has moderate growth.
It satisfies
\begin{equation}\label{stair:bounds}
 1\leq\mu_k\leq Q(k),\qquad
 \frac{\mu_k}{Q(k)}\ \text{is nonincreasing},\qquad
 Q(r_j)^{-1}\leq c_j\leq1,
\end{equation}
and
\begin{equation}\label{stair:qa-proper}
 \sum_{k\geq1}\frac1{k\mu_k}=\infty,\qquad
 M_k^{1/k}\longrightarrow\infty.
\end{equation}
Thus $\CM$ is quasianalytic, strictly contains the analytic class,
and satisfies the standing assumptions.
\end{proposition}

For the standard weight criteria, see \cite{Rainer2022Survey}.
Long constant quotient intervals are a classical way to obtain
quasianalyticity \cite{Mandelbrojt1952}; compare
\cite[\S7.4]{FurdosNenningRainerSchindl2020}, where the additional
growth properties differ.  Here their role is visible in
\[
 \sum_{k=s_j+1}^{r_{j+1}}\frac1{k\mu_k}
 \geq\frac1{Q(s_j)}\log\frac{r_{j+1}+1}{s_j+1}\longrightarrow\infty.
\]
Meanwhile $\mu_{s_j}\geq Q(s_j)/Q(r_j)\to\infty$, so the growth
intervals keep the class properly larger than the analytic class.
The proof, including moderate growth, is in
Appendix~\ref{app:stair-details}.

\begin{proposition}\label{stair:threshold}
There are constants $R>1$ and $c>0$ such that
\begin{equation}\label{stair:log-threshold}
 \Lambda_M(T)-\Lambda_M(T/R)
 \geq c\frac{T^{\beta_0}}{(\log T)^2}
 \qquad(T\text{ sufficiently large}).
\end{equation}
The weight does not satisfy $(G_{\beta_0})$.  Consequently,
\begin{equation}\label{stair:spectrum}
 \{\beta\in[0,1):M\text{ satisfies }(G_\beta)\}=[0,\beta_0).
\end{equation}
\end{proposition}

The logarithmic loss in \eqref{stair:log-threshold} is smaller than
every positive power of $T$, so every strictly subcritical condition
holds.  The failure at the endpoint is explained by the counting identity
\begin{equation}\label{stair:omega-count}
 \omega_M(T)=\sum_{\mu_k\leq T}\log\frac{T}{\mu_k}.
\end{equation}
At $T_j=\mu_{s_j}=c_jQ(s_j)$ the following plateau contributes zero.
This makes $\Lambda_M(T_j)=o(T_j^{\beta_0})$, excluding the endpoint
condition.  The counting formula is standard
\cite[Lemmas~2.4--2.5]{Schindl2024Regularization}; its verification
and the proof of the proposition are in
Appendix~\ref{app:stair-counting}.

For each fixed divisor, the polynomial degrees below will be bounded
by $E_3(Cr_j)$, with $C$ independent of $j$.
Since \eqref{stair:degree-room} holds for every fixed $C$, this one
weight accommodates all the divisors in the counterexample theorem.

\section{Failure of division at the prescribed threshold}
\label{sec:negative}\label{sec:counterexample}

Let $M$ be the weight of Section~\ref{sec:staircase}, and choose
\begin{equation}\label{neg:divisor}
 P(x,y)=y^{2m}+x^{2n},\qquad 1\leq m<n,\qquad
 q=\frac nm\geq q_*:=\frac1{1-\beta_0}=\frac\alpha{\alpha-1}.
\end{equation}
Work on $K=[-1,1]^2$.  For a positive weight $S$ and $B>0$, write
\begin{equation}\label{eq:weighted-norm}
 \|u\|_{S,B,K}
 =\sup_{\gamma\in\N^2}
 \frac{\|\partial^\gamma u\|_K}{B^{|\gamma|}|\gamma|!S_{|\gamma|}},
 \qquad \|u\|_{S,K}=\|u\|_{S,1,K}.
\end{equation}
The norm $\|u\|_{C^r(K)}$ is the maximum of the suprema of all
derivatives of total order at most $r$.  Constants may depend on
$m,n,\alpha$ and the fixed auxiliary functions, but not on $j$.

\subsection{The theorem and its geometric mechanism}
\begin{theorem}[Explicit failure of division]\label{thm:sharpness}
The weight $M$ satisfies \eqref{stair:log-threshold}.  For every pair
$(m,n)$ satisfying \eqref{neg:divisor}, there are real-valued smooth
functions $f,g$ near the origin such that
\begin{equation}\label{eq:bad-division}
 f=Pg,\qquad f\in\CM,\qquad
 g\notin\CM\text{ on any neighbourhood of }0.
\end{equation}
Moreover, there are real polynomials $p_j$, an integer $J$, and a
compact neighbourhood $K_0$ of the origin such that
\[
 P\sum_{j=J}^Np_j\longrightarrow f
 \quad\text{in the norm }\|\cdot\|_{M,K_0}.
\]
The same weight $M$ works for every such pair; the functions and $J$
may depend on the pair.
\end{theorem}

The example $\beta_0=1/2$, $P=y^2+x^4$ already exhibits the
endpoint mechanism.  Near $(a,0)$, a localized numerator can vary on
scale $a$, whereas its quotient varies on scale $a^2$ in $y$.
In general the two scales are $a$ and $a^q$.
For a cutoff $\chi$ supported near zero with $\chi(0)=1$, consider
\[
 v_a(x,y)=\frac{\chi((x-a)/a)}{y^{2m}+x^{2n}},
\]
defined as zero for $x\leq0$; its numerator is supported away from $x=0$.
Multiplication by $P$ leaves the cutoff.  On $x=a$, however, the
complex poles are $y=a^q\zeta$, $\zeta^{2m}=-1$, at distance
comparable to $a^q$ from the real axis.  This is the branch-separation
geometry of \cite{Thilliez2002Sharp}; here it is explicit in
\[
 \frac1{a^{2n}+y^{2m}}
 =\sum_{\ell\geq0}(-1)^\ell a^{-2n(\ell+1)}y^{2m\ell}.
\]
A small multiple of $v_a$ can therefore have a small numerator and
small low-order derivatives while retaining a large high-order
transverse derivative.

A smooth sum of these blocks would vanish on $\{x\leq0\}$ and be
flat at the origin.  Its nonzero numerator could not be quasianalytic.
We use the blocks only to construct polynomials: a polynomial involves
finitely many entries of the weight, so agreement of an auxiliary
weight with $M$ through its degree transfers all its estimates.
The long growth intervals provide this finite-order agreement.
Related constructions using rational blocks and large derivatives at
points approaching the origin occur in
\cite[Lemmas~3.3 and~3.11]{BelottoBierstoneKiro2022}.  The additional
requirements here are control of the numerator, the polynomial degree,
and the endpoint.

\subsection{The polynomial blocks needed for the theorem}
\label{sec:neg-blocks}
Set
\begin{equation}\label{neg:block-scales}
 a_j=e^{-r_j},\qquad z_j=(a_j,0).
\end{equation}
The weight is fixed first, the divisor second, and the blocks are
then constructed along the prescribed growth intervals.

\begin{proposition}[Polynomial blocks]\label{neg:polynomial-blocks}
For all sufficiently large $j$, there are real polynomials $p_j$ and
positive integers $k_j$ divisible by $2m$ such that
\begin{align}
 &r_j<k_j<s_j,\qquad \deg(Pp_j)<s_j,\label{neg:block-orders}\\
 &\|p_j\|_{C^{r_j}(K)}+\|Pp_j\|_{M,K}\leq2^{-j},
       \label{neg:block-smallness}\\
 &\frac{|\partial_y^{k_j}p_j(z_j)|}{k_j!M_{k_j}}
       \geq(4j)^{k_j}.\label{neg:block-large-derivative}
\end{align}
The blocks and orders are prescribed below.
\end{proposition}

No additional condition on pairs of blocks is needed.  The degree
and low-order bounds already control earlier and later summands,
respectively.  We give the construction next; the proposition is
proved in Appendix~\ref{app:subsidiary-blocks}.

\subsection{The smooth model and the cost of localization}
Introduce the three products
\begin{equation}\label{neg:aux-products}
 \mathsf A_k=\prod_{i=1}^kQ(i),\qquad
 \mathsf B_k=\prod_{i=1}^kL(i)(\log L(i))^{3/2},\qquad
 \mathsf D_k=\frac{\mathsf A_k}{\mathsf B_k},
 \qquad \mathsf A_0=\mathsf B_0=\mathsf D_0=1.
\end{equation}
The first continues the growth profile, the second supplies a smooth
cutoff, and their ratio measures the weight available after the
cutoff derivatives.  In particular,
\begin{equation}\label{neg:D-quotients}
 \frac{\mathsf D_k}{\mathsf D_{k-1}}
 =L(k)^{\alpha-1}(\log L(k))^{2\alpha-3/2}.
\end{equation}
All three quotient sequences are increasing and greater than one.

\begin{lemma}[A localized smooth model]
\label{neg:cutoffs}\label{neg:localized-estimates}
There is a real-valued $\chi\in C_c^\infty(\R)$, with
$0\leq\chi\leq1$, $\operatorname{supp}\chi\subset[-1/4,1/4]$,
and $\chi=1$ near zero, satisfying
\begin{equation}\label{neg:cutoff-bounds}
 \|\chi^{(k)}\|_\infty\leq C_\chi^{k+1}k!\mathsf B_k
 \qquad(k\geq0)
\end{equation}
for a fixed $C_\chi\geq1$.  For $0<a<1/8$, define
\begin{equation}\label{neg:localized-quotient}
 v_a(x,y)=\begin{cases}
 \displaystyle\frac{\chi((x-a)/a)}{y^{2m}+x^{2n}},&x>0,\\[4pt]
 0,&x\leq0.
 \end{cases}
\end{equation}
This function is smooth, supported in $[3a/4,5a/4]\times\R$, and,
for $k=|\gamma|$, satisfies
\begin{align}
 \|\partial^\gamma v_a\|_\infty
   &\leq C_v^{k+1}a^{-2n-qk}k!\mathsf B_k,
       \label{neg:quotient-bounds}\\
 \|\partial^\gamma(Pv_a)\|_K
   &\leq C_\chi^{k+1}a^{-k}k!\mathsf B_k.
       \label{neg:numerator-bounds}
\end{align}
Moreover $Pv_a=\chi((x-a)/a)$ everywhere, and
\begin{equation}\label{neg:exact-derivative}
 \partial_y^k v_a(a,0)=(-1)^{k/(2m)}k!a^{-2n-qk}
 \qquad(2m\mid k).
\end{equation}
\end{lemma}

The cutoff follows by infinite convolution
\cite[Lemma~3.3 and Proposition~3.11]{Rainer2022Survey}; the power
$3/2>1$ makes the convolution lengths summable.  The rescaling
\[
 x=a(1+\xi),\qquad y=a^q\eta,\qquad
 P=a^{2n}\bigl((1+\xi)^{2n}+\eta^{2m}\bigr)
\]
explains the two derivative scales.  The proof of the lemma is in
Appendix~\ref{app:localized-blocks}.

Only one auxiliary associated function is needed:
\[
 \omega_{\mathsf D}(T)=\sup_{k\geq0}\log\frac{T^k}{\mathsf D_k},
 \qquad e^{-\omega_{\mathsf D}(T)}\frac{T^k}{\mathsf D_k}\leq1
 \quad\text{for every }k.
\]
The second inequality chooses an amplitude controlling all numerator
orders at once.  Its size is given by the following estimate.

\begin{lemma}[The cost of localization]\label{neg:product-asymptotics}
As $T\to\infty$,
\begin{equation}\label{neg:D-asymptotics}
 \log\omega_{\mathsf D}(T)
 \asymp T^{1/(\alpha-1)}(\log T)^{-\nu},
 \qquad \nu=2+\frac1{2(\alpha-1)}.
\end{equation}
Consequently,
\begin{equation}\label{neg:composed-asymptotic}
 Q\bigl(\omega_{\mathsf D}(T)\bigr)
 \asymp\left(\frac{T}{\sqrt{\log T}}\right)^{q_*}.
\end{equation}
\end{lemma}
The proof is an elementary counting and inversion argument in
Appendix~\ref{app:products-cutoffs}.

\subsection{A small numerator and a large transverse derivative}
\label{sec:neg-balance}
For the estimates, continue the $j$th growth interval indefinitely:
\begin{equation}\label{neg:temporary-weight}
 \widetilde M^{(j)}_0=1,\qquad
 \widetilde M^{(j)}_k=\prod_{i=1}^k\max\{\mu_i,c_jQ(i)\}.
\end{equation}
The maximum retains the earlier quotients and continues $c_jQ(i)$
after $s_j$.
\begin{lemma}[Temporary weights]\label{neg:temporary-properties}
The sequence $\widetilde M^{(j)}$ is increasing and logarithmically
convex, its associated class is nonquasianalytic, and
\begin{equation}\label{neg:temporary-bounds}
 \widetilde M^{(j)}_k=M_k\quad(k\leq s_j),\qquad
 \widetilde M^{(j)}_k\geq\max\{1,c_j^k\mathsf A_k\}\quad(k\geq0).
\end{equation}
\end{lemma}
Its proof is in Appendix~\ref{app:temporary-weights}.

Choose the amplitude, and then the derivative order, by
\begin{equation}\label{neg:amplitudes}
 \eta_j=2^{-j-10},\qquad
 \tau_j=\frac{C_\chi}{c_ja_j},\qquad
 w_j=\frac{\eta_j}{C_\chi}e^{-\omega_{\mathsf D}(\tau_j)},\qquad
 u_j=w_jv_{a_j},
\end{equation}
\begin{equation}\label{neg:testing-order}
 k_j=2m\left\lceil\frac{\log(1/w_j)}{2m}\right\rceil.
\end{equation}
The amplitude suppresses all numerator orders; the choice of $k_j$
ensures $w_j^{1/k_j}\geq e^{-1}$.  Thus the amplitude costs only a
bounded factor in the $k_j$th root of the normalized derivative.

\begin{lemma}[The weighted smooth blocks]
\label{neg:small-estimates}\label{neg:large-derivative}
For every $j$,
\begin{equation}\label{neg:small-numerator}
 \|Pu_j\|_{\widetilde M^{(j)},K}\leq\eta_j.
\end{equation}
For all sufficiently large $j$,
\begin{equation}\label{neg:small-low-orders}
 \|u_j\|_{C^{r_j}(K)}\leq\eta_j,
\end{equation}
and
\begin{equation}\label{neg:large-derivative-bound}
 r_j<k_j<s_j,\qquad
 \frac{|\partial_y^{k_j}u_j(z_j)|}{k_j!M_{k_j}}\geq(8j)^{k_j}.
\end{equation}
More precisely, for a constant $C$ independent of $j$,
\begin{equation}\label{neg:endpoint-gap}
 \left(\frac{|\partial_y^{k_j}u_j(z_j)|}{k_j!M_{k_j}}\right)^{1/k_j}
 \geq C^{-1}e^{(q-q_*)r_j}
       \left(\frac{\sqrt{r_j}}{Q(r_j)}\right)^{q_*}.
\end{equation}
\end{lemma}

\begin{remark}[The endpoint gain]\label{neg:critical-remark}
Above the threshold, the exponential in \eqref{neg:endpoint-gap}
grows.  At $q=q_*$ it equals one, but the remaining factor still
grows faster than $j$, since $Q(r)$ grows only through logarithms
and $j\leq\log r_j$.  The power $3/2$ in the cutoff product lies
between $1$ and $2$: the first inequality permits a cutoff, and the
second leaves the factor $r_j^{q_*/2}$ responsible for the endpoint.
\end{remark}
The proof uses only $M_k\leq Q(k)^k$, the exact derivative, and
\eqref{neg:composed-asymptotic}; see
Appendix~\ref{app:weighted-blocks}.

\subsection{Polynomial replacement and finite-order agreement}
The smooth blocks now have all the required estimates.  To obtain
quasianalytic numerators, replace each by a polynomial while retaining
its low derivatives and its selected large derivative.  The essential
observation is
\begin{equation}\label{neg:finite-degree-transfer}
 \deg(Pp)<s_j\quad\Longrightarrow\quad
 \|Pp\|_{M,K}=\|Pp\|_{\widetilde M^{(j)},K}.
\end{equation}
The weights agree at every possibly nonzero derivative order of $Pp$.
No approximation of a localized smooth block in the final
quasianalytic norm is asserted.

\begin{lemma}[Quantitative polynomial replacement]
\label{neg:gaussian-approximation}
There are prescribed real Gaussian--Taylor polynomials $p_j$ and
integers $d_j$ such that, for all sufficiently large $j$,
\begin{equation}\label{neg:approximation-bound}
 \|p_j-u_j\|_{\widetilde M^{(j)},K}\leq\delta_j:=e^{-r_j^2},\qquad
 \deg p_j\leq2d_j,\qquad d_j\leq E_3(Cr_j),
\end{equation}
where $C$ is independent of $j$.  In particular,
\begin{equation}\label{neg:stage-approximation}
 \deg(Pp_j)\leq2d_j+2n<s_j.
\end{equation}
The Gaussian widths and truncation degrees are given explicitly in
\eqref{neg:approximation-prescription}.
\end{lemma}

Gaussian smoothing and Taylor truncation give the polynomials;
compare \cite{Plesniak1994} for related 
approximation.  The proof in Appendix~\ref{app:negative-approx}
controls all derivative orders in one weighted norm as well as the
degree.  Global boundedness, rather than compact support, suffices
for the Gaussian convolution to be entire.

We also use the elementary multiplication estimate, for increasing $S$,
\begin{equation}\label{neg:polynomial-multiplication}
 \|Fu\|_{S,K}\leq C_F\|u\|_{S,K},\qquad
 C_F=\sum_\delta\frac{\|\partial^\delta F\|_K}{\delta!}
 \quad(F\text{ a fixed polynomial}).
\end{equation}
Together with \eqref{neg:finite-degree-transfer}, the approximation
bound gives
\begin{equation}\label{neg:polynomial-transfer}
 \|p_j-u_j\|_{C^{r_j}(K)}\leq\eta_j,\qquad
 \|Pp_j\|_{M,K}\leq2\eta_j
\end{equation}
for large $j$.  Since $k_j<s_j$, the large derivative is preserved
as well.  These transfers prove Proposition~\ref{neg:polynomial-blocks};
the details are in Appendix~\ref{app:subsidiary-blocks}.

\subsection{Summation and noncancellation}
For a sufficiently large integer $J$, set
\begin{equation}\label{neg:final-functions}
 g=\sum_{j\geq J}p_j,\qquad f=Pg.
\end{equation}
The two smallness bounds give smooth convergence of $g$ and convergence
of the numerator series in one fixed $M$-norm.  At order $k_j$,
earlier polynomials vanish because their degrees are smaller than
$k_j$; later polynomials are small because $k_j<r_i$ for $i>j$.
Thus the large derivative survives.  Finally, $P>0$ off the origin
makes the smooth quotient unique.  The complete proof of
Theorem~\ref{thm:sharpness}, including fixed-norm convergence, is in
Appendix~\ref{app:negative-summation}.

\section{Division by quasianalytic germs and Noetherianity}
\label{sec:algebra}
The goal of this section is to prove the positive Noetherianity result stated in
Theorem~\ref{thm:plane}. The main step is to extend the formal division property for analytic germs, established in Corollary~\ref{cor:formal}, to arbitrary quasianalytic germs in two variables. Throughout this section, fix:
\begin{equation}\label{eq:HypSection}
\text{a weight sequence $M$ satisfying $(G_\beta)$ for every
$0\leq\beta<1$,}
\end{equation}
and we write
\[
A:=A_{M,2},\qquad B:=B_2=\mathbb R[[x,y]],
\qquad
\mathfrak m=(x,y)A,\quad \mathfrak n=(x,y)B.
\]
Note that reciprocals of nonvanishing germs belong to the class, so $A$ is local
with maximal ideal $\mathfrak m$. We identify $A$ with its image under the Taylor
expansion map. 

Recall that the ring of formal power series $B$ is a regular local ring of dimension $2$; in particular it is a unique factorization domain. As a consequence, the notion of an element being \emph{reduced} (also called squarefree) is well-defined over $B$. Over $A$, the notion of being reduced is not (yet) available, and it will be useful to consider:

\begin{definition}\label{def:squarefree}
    A nonzero germ $\varphi \in A$ is \emph{formally reduced} if its Taylor series $\hat{\varphi}$ is reduced in $B$.  Units are included in this convention.
\end{definition}

Now, recall that Corollary~\ref{cor:formal} gives a formal division \eqref{eq:formal-division} property (c.f. Remark \ref{rem:ClosedIdealsDef}):
\begin{equation}\label{eq:formal-division}
\varphi A=A\cap \widehat\varphi B
\qquad
\text{for every }0\ne \varphi\in\mathbb R\{x,y\}.
\tag{FD}
\end{equation}
We will show that the same identity holds when the analytic germ is replaced by an arbitrary nonzero germ
$\varphi \in A$, that is, our goal is to establish:

\begin{proposition}\label{prop:all-principal-contraction}
Every nonzero $\varphi\in A$ satisfies
\begin{equation}
 \widehat\varphi B\cap A=\varphi A.
 \label{eq:all-principal-contraction}
\end{equation}
Moreover, $\varphi$ is a unit times a finite product of powers of formally reduced elements in $A$.
\end{proposition}

By combining this result with \cite[Proposition~3]{Valabrega1973}, we deduce:

\begin{proof}[Proof of Theorem \ref{thm:plane}]
Our goal is to show that we can apply the criteria from
\cite[Proposition~3]{Valabrega1973}. To this end, let $C=A_{M,1}$, which is a discrete
valuation ring (in fact, every nonzero one-variable germ $c\in C$ is of the form $x^r$ times a unit: quasianalyticity makes the order finite and repeated coordinate division gives the factorization). Stability under differentiation gives an inclusion
\(
A\subset C[[y]]\subset B.
\)
Now, we claim that:
\[
A=\operatorname{Frac}(A)\cap C[[y]].
\]
Note that the Proposition follows from this claim, combined with \cite[Proposition~3]{Valabrega1973}.

To prove the claim, note that it is clear that $A\subset \operatorname{Frac}(A)\cap C[[y]]$, so consider $u \in \operatorname{Frac}(A)\cap C[[y]] $. We know that $u = f/g$, where $f$ and $g\in A$. Moreover, since $u$ belongs to $C[[y]] \subset B$, we have that $\widehat f = \widehat g\,H$, where $H=u\in B$. From the fact that principal ideals are closed in $A$, we get that $f =g a$ for some $a\in A$ and, therefore, $u = a$ as needed.
\end{proof}

\subsection{Proof of Proposition \ref{prop:all-principal-contraction}: formally reduced germs}
The goal of this subsection is to establish the following version of Proposition \ref{prop:all-principal-contraction}:

\begin{lemma}\label{lem:reduced-division-lifting}
Let $\varphi \in A$ be nonzero and
formally reduced. Then
\begin{equation}
 \widehat \varphi B\cap A=\varphi A.
 \label{eq:reduced-contraction}
\end{equation}
Moreover, every formal factor of $\widehat \varphi$
is associated to the Taylor series of an element in $A$.
\end{lemma}

The proof of this Lemma is based on a key technical result, which is essentially a $\CM$-version of a well-known ``finite determinacy result", see for instance \cite[Chapter VIII, §5, Proposition 5.1.2]{Tougeron1972}. The extra $\CM$ regularity is obtained thanks to \cite{RainerSchindl2016}:

\begin{lemma}[Finite determinacy]\label{lem:reduced-normal-form}
Let $\varphi\in\mathfrak m$ be formally reduced. Then, for every sufficiently large integer $N$, there exists a local $\CM$-diffeomorphism $\Psi_N$, tangent to the identity at the origin, such that
\[
\varphi=P_N\circ\Psi_N,
\]
where $P_N$ is the Taylor polynomial of $\varphi$ of degree $N$. In particular, $P_N$ is reduced as an analytic germ.
\end{lemma}

The proof of this result can be found in Appendix \ref{app:reduced-normal-form}. We are ready to prove Lemma \ref{lem:reduced-division-lifting}:

\begin{proof}[Proof of Lemma \ref{lem:reduced-division-lifting}]
By Lemma \ref{lem:reduced-normal-form}, write $\varphi=P\circ\Psi$, where $P$ is a reduced polynomial and $\Psi$ is a local $\CM$-diffeomorphism. The first assertion follows
from property \eqref{eq:formal-division}, after composition with $\Psi^{-1}$. The second follows from the first argument, combined with the fact that formal factors of analytic plane germs are convergent up to units, by Newton--Puiseux theory
\cite[\S4]{Teissier2007}.
\end{proof}

\subsection{Proof of Proposition \ref{prop:all-principal-contraction}: the general case}

The main idea is to reduce the general case to a recursive use of Lemma \ref{lem:reduced-division-lifting}. This reduction is done via a standard formal Bertini-type argument:

\begin{lemma}\label{lem:derivative-pencil}
Let $0\ne F\in\mathfrak n\subset B$, and write
\[
F=U\prod_{j=1}^s F_j^{m_j},
\qquad U \text{ is a unit},
\]
where the $F_j$ are pairwise nonassociated irreducible elements of $B$.
Set
\[
\ell=\max_j m_j,
\qquad
Q=\prod_{m_j=\ell}F_j.
\]
Then there are real constants $c_\alpha$, $|\alpha|<\ell$,
such that
\[
\sum_{|\alpha|<\ell}c_\alpha D^\alpha F
\]
is reduced and divisible by $Q$.
\end{lemma}

\begin{proof}
In characteristic zero, differentiation gives
\[
\gcd\{D^\alpha F:|\alpha|<\ell\}\sim Q.
\]
Indeed, writing $F=H^eV$ for an irreducible
factor $H\nmid V$ and $1\leq e\leq \ell$, choose a coordinate derivative
$\partial$ such that $H\nmid \partial H$. By induction,
\[
\partial^jF\equiv
\frac{e!}{(e-j)!}H^{e-j}(\partial H)^jV
\pmod{H^{e-j+1}},
\qquad 0\le j\le e.
\]
Thus derivatives of order at most $r < e\leq \ell$ have
greatest common divisor of multiplicity
$\max(e-r,0)$ along $H$, proving the claim.

Next, writing $D^\alpha F=QG_\alpha$, we therefore
have $\gcd\{G_\alpha; |\alpha| < \ell\}=1$. We now claim that a general real linear
combination
\[
G=\sum_{|\alpha|<\ell}c_\alpha G_\alpha
\]
is reduced and coprime with $Q$. Indeed, if the span of the $G_\alpha$ contains a unit,
there is nothing to prove. Otherwise, choose two
coprime elements $G_1,G_2$ in this span.
For some coordinate derivative $\partial$,
\[
S=G_2\partial G_1-G_1 \partial G_2\neq 0.
\]
Given $t\in \mathbb{R}$, every repeated irreducible factor of
$G_1+tG_2$ divides $S$. Since $G_1$ and $G_2$
are coprime, each irreducible factor of $SQ$
divides $G_1+tG_2$ for at most one
$t\in\mathbb R$. Hence, for generic $t$,
$G_1+tG_2$ is reduced and coprime to $Q$.

Thus $QG$ is reduced and gives the
required combination of derivatives.
\end{proof}

We are now ready to prove:

\begin{proof}[Proof of Proposition \ref{prop:all-principal-contraction}]
Suppose $\varphi$ is a nonunit and apply
Lemma~\ref{lem:derivative-pencil} to
$F=\widehat\varphi$, with the notation used there.
The same real coefficients applied to the actual
derivatives give a germ
\[
b=\sum_{|\alpha|<\ell}c_\alpha D^\alpha\varphi\in A
\]
whose Taylor series is reduced and divisible by $Q$. By Lemma~\ref{lem:reduced-division-lifting},
there exists a divisor $\rho\in A$ of $b$ such that
\[
\widehat\rho B=QB.
\]
In particular, $\rho$ is formally reduced and
$\widehat\rho^\ell$ divides $\widehat\varphi$.
Repeated application of
\eqref{eq:reduced-contraction} gives
\[
\varphi=\rho^\ell\varphi_1,
\qquad \varphi_1\in A.
\]
The maximal formal multiplicity of $\varphi_1$
is strictly less than $\ell$, unless $\varphi_1$
is a unit. Induction gives the asserted
factorization of $\varphi$. Finally, applying
\eqref{eq:reduced-contraction} successively
to the factors of $\varphi$ gives
\(
\widehat\varphi B\cap A=\varphi A.
\)
\end{proof}

 \subsection{Smooth definable quotients} Let $\R_M$ be the expansion of the real field generated by
 restricted $\CM$ functions.  It is o-minimal and polynomially
 bounded~\cite{RolinSpeisseggerWilkie2003}.  Recall the ring
 \[
  \mathcal D_d(M)=
  \{u\in C^\infty_0(\R^d):u\text{ is definable near }0
           \text{ in }\R_M\}.
 \]
 It is a quasianalytic enlargement of $A_{M,d}$;
 see~\cite{Miller1995} and
 \cite[\S2.2]{BelottoBiborskiBierstone2017}.

 \begin{proposition}\label{prop:definable-quotient}
 If $f,\varphi\in A_{M,d}$, $\varphi\ne0$, and $f=\varphi g$ for a
 smooth germ $g$, then $g\in\mathcal D_d(M)$.
 In particular, each quotient in Theorem~\ref{thm:sharpness}
 belongs to $\mathcal D_2(M)\setminus A_{M,2}$.
 \end{proposition}

 \begin{proof}
 Choose a small definable box $U$ on which $f,\varphi,g$ are
 defined, and on which $f$ and $\varphi$ are definable.
 The graph of $f/\varphi$ over $U\setminus\{\varphi=0\}$ is
 definable.  This domain is dense in $U$, by quasianalyticity.
 Continuity of $g$ shows that its graph is the relative closure of
 that graph in $U\times\R$, and is therefore definable.
 \end{proof}

 More generally, Nowak's criterion for divisibility
 \cite{Nowak2013} gives formal principal contraction in the full
 smooth definable ring $\mathcal D_d(M)$, in every dimension.
 The counterexamples above distinguish division within a fixed
 Denjoy--Carleman class from division within this definable
 quasianalytic enlargement.  For their scalar equations
 $P(x,y)u-f(x,y)=0$, a formal solution is realized by a smooth
 definable solution, while no solution belongs to the original
 class.

\subsection{Nonclosed ideal in the Silva topology}

Finally, we provide the small argument necessary to justify Remark \ref{rem:Silva}:

\begin{corollary}\label{cor:closednessBoth}
For the weight $M$ and the polynomial $P$ in
Theorem~\ref{thm:sharpness}, the principal ideal
$P A_{M,2}$ is not closed in the Silva topology.
\end{corollary}

\begin{proof}
Let $f=Pg$ be the counterexample from
Theorem~\ref{thm:sharpness}, with $f\in A_{M,2}$
and $g\notin A_{M,2}$. Then $f\notin P A_{M,2}$. The polynomial multiples constructed in that
theorem converge to $f$ in a fixed $\CM$-norm,
hence in the Silva topology. Thus $P A_{M,2}$
is not Silva-closed.
\end{proof}

\appendix

\section{Holomorphic approximation and auxiliary estimates}
\label{app:holomorphic}

We prove the approximation statements used in
Section~\ref{sec:positive}, and collect the auxiliary estimates for
central smallness, continuation, Cauchy division, and removal of a
small factor. Throughout this appendix, the standing
assumptions on $M$ hold, $M_k^{1/k}\to\infty$, and $M$ has moderate
growth. We use the polydiscs $D_{\C}^{j}(p,r)$, their real parts
$D_{\R}^{j}(p,r)$, and the tubes $T_r(K)$ introduced in that section.
The only dimensions needed are $j=1,2$. As in Section~\ref{sec:positive},
we abbreviate $D_r(p)=D_{\C}^{2}(p,r)$ and
$D_r^{\R}(p)=D_{\R}^{2}(p,r)$, omitting $p$ when $p=0$.

\subsection{From almost analytic extensions to bounded approximants}

Dynkin's almost analytic extension theorem \cite{Dynkin1980}
(see also \cite[\S1.2, Theorem~4]{Dynkin1993}), in the
several-variable formulation of
\cite[Theorem~1.2]{FurdosNenningRainerSchindl2020}, characterizes
$\CM$ by extensions whose $\overline\partial$-derivatives decay at
the rate $h_M$. We use the extension theorem as an input; the short
argument below converts its conclusion into holomorphic functions
on shrinking tubes. For approximation characterizations on
intervals, see also \cite[Theorem~5.3]{NenningRainerSchindl2021}.

\begin{proof}[Proof of Proposition~\ref{prop:hol-approx}]
Apply the extension theorem to $u$ on the ball $B\Subset U$, with
full derivative weight $\mathcal M_k=k!M_k$. It gives a bounded
$C^1$ extension $V$ such that $V|_B=u$ and
\begin{equation}\label{eq:dynkin-local}
 |\overline\partial V(z)|\leq C h_M(A|\operatorname{Im}z|)
\end{equation}
when the real projection of $z$ belongs to $B$. The hypotheses of
the cited theorem are logarithmic convexity and root growth of the
normalized weight $\mathcal M_k/k!=M_k$, together with boundedness
of $(\mathcal M_{k+1}/\mathcal M_k)^{1/k}$. The last condition follows
from differentiation stability, since
\[
 \left(\frac{\mathcal M_{k+1}}{\mathcal M_k}\right)^{1/k}
 \leq (k+1)^{1/k}D^{(k+1)/k}\qquad(k\geq1).
\]
Complex-valued $u$ are treated componentwise.

Choose a real rectangle $K_1$ with
$K\subset\operatorname{int}K_1\subset K_1\subset B$. For small $s$,
let $\mathcal R_s=R_{1,s}\times R_{2,s}$ be the product of complex
rectangles with real projection $K_1$ and imaginary half-widths
$(T+2)s$. The tube $T_{Ts}(K)$ lies in $\mathcal R_s$, at distance
at least $c_Ts$ from each coordinate boundary. On $\mathcal R_s$,
\eqref{eq:dynkin-local} is bounded by $C h_M(A_Ts)$.

Let $C_j$ denote the boundary Cauchy operator and $S_j$ the solid
Cauchy--Green operator in the $j$th variable, with signs chosen so
that $V-C_jV=S_j\partial_{\bar z_j}V$. Set $u_s=C_1C_2V$. The
boundary integral defines a jointly holomorphic function on
$\mathcal R_s$. Applying Cauchy--Green in each variable gives
\[
 V-u_s=S_1\partial_{\bar z_1}V
          +C_1S_2\partial_{\bar z_2}V.
\]
The area integral of $|\zeta-z|^{-1}$ over each rectangle is bounded
uniformly in $s$, so the operators $S_j$ have uniformly bounded
supremum norms. On $T_{Ts}(K)$ the boundary integral gives
$\|C_1\|\leq C_Ts^{-1}$. These estimates also apply when the other
variable ranges over a boundary of its rectangle. Thus
\[
 \sup_{T_{Ts}(K)}|V-u_s|\leq C_Ts^{-1}h_M(A_Ts).
\]
Lemma~\ref{lem:weight-loss} absorbs $s^{-1}$. Restricting to $K$
gives the required approximation error, and boundedness of $V$
gives the uniform bound for $u_s$.
\end{proof}

\subsection{Elementary continuation and the branch geometry}

The first lemma is the two-constants theorem applied to a slit
domain. Its role is to turn smallness on a real square into
smallness on a smaller complex polydisc. The corresponding
one-variable estimates on shrinking ellipses are
\cite[Lemma~3.2.4]{Thilliez2020Powers} and
\cite[Lemma~5.2]{NenningRainerSchindl2021}. The following is the
fixed-shape, one- and two-variable version needed here.

\begin{lemma}[Real-to-complex smallness]\label{lem:real-complex}
For $j=1,2$ there is $\theta_j\in(0,1)$ with the following property.
If $H$ is holomorphic on $D_{\C}^{j}(0,2)$,
$|H|\leq B$ there, and $|H|\leq\varepsilon$ on $[-1,1]^j$, where
$0<\varepsilon\leq B$, then
\[
 \sup_{D_{\C}^{j}(0,1)}|H|
 \leq B^{1-\theta_j}\varepsilon^{\theta_j}.
\]
Translations, dilations, and other fixed choices of nested radii
give the same type of estimate, with an exponent independent of
the centre and scale.
\end{lemma}

\begin{proof}
For $j=1$, let $\omega$ be the harmonic measure of the slit
$[-1,1]$ in $D_{\C}^{1}(0,2)\setminus[-1,1]$. The two-constants
theorem for the subharmonic function $\log(|H|/B)$ gives
\[
 \log(|H(z)|/B)\leq\omega(z)\log(\varepsilon/B).
\]
The harmonic measure is positive in the slit domain and extends
with value one to the slit. Its endpoints are regular boundary
points, as a square-root coordinate gives a local barrier.
Consequently, $\omega$ has a positive lower bound $\theta_1$ on
the closed unit disc outside the slit. This proves the estimate;
on the slit it follows from the hypothesis. Zeros of $H$ cause no
difficulty in the subharmonic argument.

For $j=2$, first fix the second variable real and apply the
one-variable estimate in the first variable. Then fix the first
variable in the complex unit disc and apply it in the second
variable. This gives $\theta_2=\theta_1^2$. The same
harmonic-measure argument with fixed radii, followed by translation
and dilation, proves the last assertion.
\end{proof}

The second lemma is an iteration of Hadamard's three-circles
theorem. For the quantitative continuation viewpoint, see also
\cite[Theorem~3.1 and \S7]{Trefethen2020Continuation}. In
Section~\ref{sec:positive}, the number of discs is determined by the
geometry of a Puiseux branch.

\begin{lemma}[Propagation along discs]\label{lem:disc-chain}
Let $z_0,\ldots,z_N\in\C$ satisfy $|z_{i+1}-z_i|\leq r$.
Suppose $H$ is holomorphic and bounded by $B$ on
$\bigcup_{i=0}^{N}D_{\C}^{1}(z_i,4r)$, and is bounded by
$0<\varepsilon\leq B$ on $D_{\C}^{1}(z_0,r)$. Then
\[
 \sup_{D_{\C}^{1}(z_N,r)}|H|
 \leq B^{1-2^{-N}}\varepsilon^{2^{-N}}.
\]
\end{lemma}

\begin{proof}
At the radii $r,2r,4r$, the three-circles theorem gives
\[
 \sup_{D_{\C}^{1}(z_i,2r)}|H|
 \leq B^{1/2}
       \left(\sup_{D_{\C}^{1}(z_i,r)}|H|\right)^{1/2}.
\]
Since $D_{\C}^{1}(z_{i+1},r)\subset D_{\C}^{1}(z_i,2r)$,
iteration proves the assertion.
\end{proof}

\paragraph{Geometric details for Lemma~\ref{pos:branch-smallness}.}
\label{app:branch-geometric-details}
Using the sectorial branch estimates in \eqref{pos:geometry-used},
cf.\ \cite[\S\S1.3--1.4]{Thilliez2002Sharp}, we record the choices
that ensure that the graph discs used in the main proof remain in
the approximation tube. All these choices are
made independently of $L$. Choose $\eta>0$ so small that discs of
radius $4\eta s$ centred at real points of absolute value at least
$s$ lie in the relevant Puiseux sector. By
\eqref{pos:geometry-used}, their graphs lie in an $O(s)$-tube of the
real plane whenever the central root has height $O(s)$. Choose
$K_0$ large enough to contain the graph discs centred at $\pm s$,
all graphs with $\sigma\le|a|\le s$ over $|z-a|<\eta\sigma$, and
the central contours in \eqref{pos:central-integral}. These contours
have radius $R_0\sigma$ for a constant $R_0$ depending only on the
root bound in \eqref{pos:geometry-used}; since $\sigma\le s$, this
places no dependence on $L$. Choose $T$ large enough to contain all
these constructions and $D_{2K_0s}$. Shrink $r$ so that their real
projections remain in the interior of $K$.

If $\sigma\le|a|\le s$, the graph lies in $D_{K_0s}$, and
Lemma~\ref{pos:central} gives the required bound with $h_M(A_{\mathrm{cen}}s)$.
Once $L$ is fixed, this has the form $h_M(A\sigma)$.

If $|a|\ge s$ and the branch is real on its half-axis, put
$H(z)=F_s(z,\lambda_\nu(z))$. This is holomorphic and uniformly
bounded on $|z-a|<4\eta s$. On the real interval $|t-a|<\eta s$,
\eqref{pos:identity} gives $f(t,\lambda_\nu(t))=0$; hence
$|H(t)|\le C h_M(A_1s)$. Lemmas~\ref{lem:real-complex} and
\ref{lem:weight-loss} give the same form of smallness on
$|z-a|<\eta s$, after a fixed decrease of $\eta$ if necessary.

For the nonreal branch treated in the main proof, the real centres
between $\operatorname{sgn}(a)s$ and $a$ have absolute value at
least $s$. By \eqref{pos:height-monotone}, their central heights
are at most $C\sigma$. The derivative bound in
\eqref{pos:geometry-used} adds at most $C\eta s$ to the imaginary
height on each disc. Thus the graphs of all discs of radius
$4\eta s$ lie in $T_{Ts}(K)$, and the first graph disc of radius
$\eta s$ lies in $D_{K_0s}$. These facts justify the uniform
boundedness and the initial smallness used in the disc-chain
argument.

\subsection{Recovering derivative bounds from local approximants}

Proposition~\ref{prop:approx-characterization} is a local,
two-variable variant of
\cite[Proposition~3.3.2]{Thilliez2020Powers} and
\cite[Theorem~5.3(ii)]{NenningRainerSchindl2021}. The local criterion
allows polynomial growth of the complex
supremum, and requires no compatibility between approximants at
different centres. We compare successive scales at the same centre.
The following proof also explains why the resulting derivative
limits agree with the derivatives of $u$, despite the shrinking
domains.

\begin{proof}[Proof of Proposition~\ref{prop:approx-characterization}]
Choose $0<\sigma_*<\min\{\sigma_0,1\}$ and put
$\sigma_i=\sigma_*2^{-i}$. Fix $p\in K$ temporarily and write
$U_i=U_{p,\sigma_i}$. On the common real square of radius
$\kappa\sigma_i/2$, the hypothesis gives
\[
 |U_{i+1}-U_i|\leq 2C h_M(B\sigma_i).
\]
On the common complex polydisc of radius $\kappa\sigma_i/2$, the
difference is bounded by $C_0\sigma_i^{-a}$. Restrict its real
smallness estimate to the square of radius $\kappa\sigma_i/4$ and
apply Lemma~\ref{lem:real-complex} with outer radius
$\kappa\sigma_i/2$ and inner radius $\kappa\sigma_i/4$.
Absorb the fixed fractional power and polynomial
factor by Lemma~\ref{lem:weight-loss}. We obtain constants
$c_0,B_0,C_1>0$, independent of $i$ and $p$, such that
\begin{equation}\label{eq:hol-approx-difference}
 \sup_{D_{\C}^{2}(p,c_0\sigma_i)}|U_{i+1}-U_i|
 \leq C_1h_M(B_0\sigma_i).
\end{equation}
Reducing $\sigma_*$ if necessary makes the smallness hypothesis
hold at every scale. If $n=|\alpha|$, Cauchy's inequality gives
\begin{equation}\label{eq:hol-approx-derivative-difference}
 |\partial^\alpha(U_{i+1}-U_i)(p)|
 \leq C_1\alpha!(c_0\sigma_i)^{-n}h_M(B_0\sigma_i).
\end{equation}
For each fixed $\alpha$, this is summable in $i$.

To identify the derivative limit, let $T_m u$ be the Taylor
polynomial of $u$ at $p$, of degree $m$. Smoothness and the fact
that $h_M(t)$ is smaller than every power of $t$ imply
\[
 \sup_{D_{\R}^{2}(p,\kappa\sigma_i)}|U_i-T_m u|
 \leq C_m\sigma_i^{m+1}.
\]
On the corresponding complex polydisc, the supremum is at most
$C_m\sigma_i^{-a}$. Lemma~\ref{lem:real-complex}, followed by
Cauchy's inequality on a fixed smaller polydisc, therefore gives
\[
 |\partial^\alpha U_i(p)-\partial^\alpha u(p)|
 \leq C_{m,\alpha}
       \sigma_i^{\theta(m+1)-(1-\theta)a-n}
 \qquad(m\geq n)
\]
for a fixed $\theta\in(0,1)$. Choose $m$ large enough to make the
exponent positive. It follows that
$\partial^\alpha U_i(p)\to\partial^\alpha u(p)$. The constants in
this identification may depend on $p$; they do not enter the
uniform bound that follows.

Telescoping \eqref{eq:hol-approx-derivative-difference} gives
\[
 |\partial^\alpha u(p)|
 \leq C_2C_3^n\alpha!
      +C_1c_0^{-n}\alpha!
       \sum_{i\geq0}\sigma_i^{-n}h_M(B_0\sigma_i).
\]
The first term is Cauchy's estimate for $U_0$, on its polydisc of
fixed radius. By the definition of $h_M$,
\[
 \sum_{i\geq0}\sigma_i^{-n}h_M(B_0\sigma_i)
 \leq B_0^{n+1}M_{n+1}\sum_{i\geq0}\sigma_i
 \leq C_4C_5^nM_n.
\]
Here $\sum_i\sigma_i=2\sigma_*$, and the last inequality is
differentiation stability. Since $M_n\geq1$, the analytic first
term is absorbed as well. All constants in the resulting bound
are independent of $p\in K$, which proves the proposition.
\end{proof}

\subsection{Formal quotient coefficients and optimal truncation}

The next lemma concerns formal series only. It supplies the
coefficient bound used to truncate the quotient in the central
estimate. The general closure result for formal division is
\cite[Theorem~1]{Thilliez2002Formal}; the elementary proof below
uses the adapted form of $W$ and differentiation stability.

\begin{lemma}\label{pos:formal-coefficients}
Let $W$ be an adapted Weierstrass polynomial as in
\eqref{eq:adapted-weierstrass}, and let $f\in\CM$ near the origin.
If $\widehat f=W\widehat g$, there are $A,B>0$ such that
\begin{equation}\label{pos:formal-bound}
 \widehat g=\sum_{\alpha\in\N^2}g_\alpha X^\alpha,
 \qquad |g_\alpha|\leq A B^{|\alpha|}M_{|\alpha|}
\end{equation}
\end{lemma}

\begin{proof}
Write $\widehat f=\sum_{j\geq m}F_j$ and
$\widehat g=\sum_{n\geq0}g_n$ in homogeneous components.
The order of $W$ is $m$, and its initial homogeneous part satisfies
$W_m(0,1)=1$. Hence
\[
 W(tu,t)=t^mA(t,u),\qquad A(0,0)=1,
\]
where $A$ is analytic. Choose a bidisc on which $A$ is invertible.
For a smaller closed $u$-disc $|u|\leq\varepsilon<1$, the expansion
\[
 A(t,u)^{-1}=\sum_{\ell\geq0}a_\ell(u)t^\ell
\]
satisfies $|a_\ell(u)|\leq C D_0^\ell$. Comparing coefficients of
$t$ in the formal identity gives
\[
 g_n(u,1)=\sum_{j=m}^{n+m}F_j(u,1)a_{n+m-j}(u).
\]
The derivative bounds for $f$ and monotonicity of $M$ imply
\[
 \sup_{|u|\leq\varepsilon}|g_n(u,1)|
 \leq C E^n M_{n+m}.
\]
The numbers of coefficients and summands contribute only a fixed
exponential factor in $n$. Since $g_n(u,1)$ has degree at most $n$,
Cauchy's inequality bounds each of its coefficients by
$C(E/\varepsilon)^nM_{n+m}$. Finally, differentiation stability gives
\[
 M_{n+m}\leq D^{mn+m(m+1)/2}M_n.
\]
Absorbing the fixed constants proves \eqref{pos:formal-bound}.
\end{proof}

\paragraph{Optimal truncation and central smallness.}
\label{app:central-proof}
The Taylor-remainder decomposition below is an analogue of
\cite[Lemma~4, inequalities~(20)--(21)]{Thilliez2002Formal}.
There the divisor is polynomial. Here its analytic coefficients
produce a convergent tail; optimal truncation and
Lemma~\ref{lem:real-complex} then give the complex estimate used in
Lemma~\ref{pos:central}.

\begin{proof}[Proof of Lemma~\ref{pos:central}]
Choose a large fixed $C_*$, to be increased below, and let $N=N(s)$
attain the infimum
\[
 h_M(C_*s)=(C_*s)^NM_N.
\]
Then $N(s)\to\infty$ as $s\downarrow0$. Logarithmic convexity gives
\begin{equation}\label{pos:optimum-intermediate}
 (C_*s)^nM_n\le1\qquad(0\le n\le N).
\end{equation}
Set $Q_s=T_{N-1}\widehat g$, where $T_{N-1}$ denotes truncation in
total degree. The formal coefficient bound gives, on $D_{2K_0s}$,
\[
 |Q_s|\le A_0\sum_{n<N}(n+1)(2BK_0s)^nM_n
 \le A_0\sum_{n<N}(n+1)(2BK_0/C_*)^n.
\]
This is uniformly bounded when $C_*>4BK_0$.

Put $V_s=WQ_s$ and $\rho=K_0s$. Formal divisibility gives
$T_{N-1}f=T_{N-1}V_s$. Taylor's formula on real segments gives
\begin{equation}\label{pos:taylor-remainder}
 \|f-T_{N-1}f\|_{D_\rho^{\R}}
 \le C(C_0\rho)^NM_N.
\end{equation}
For the analytic tail of $V_s$, write $W=\sum_{j\ge0}W_j$ in
homogeneous components. Its coefficients are bounded by $CD_1^j$,
so, for small $\rho$,
\begin{align*}
 \|V_s-T_{N-1}V_s\|_{D_\rho}
 &\le C\sum_{n<N}(n+1)(B\rho)^nM_n
       \sum_{j\ge N-n}(j+1)(D_1\rho)^j\\
 &\le C\rho^NM_N
       \sum_{n<N}(n+1)(N-n+1)B^nD_1^{N-n}\\
 &\le C(C_1\rho)^NM_N.
\end{align*}
Here $M_n\le M_N$, and the last finite sum is bounded by a fixed
exponential in $N$.

Increase $C_*$ to dominate $C_0K_0$, $C_1K_0$ and the constant $A$
in \eqref{pos:approximants}. Combining the two remainders with
$F_s-f$ gives
\[
 \|F_s-V_s\|_{D_{K_0s}^{\R}}\le C h_M(C_*s).
\]
Both $F_s$ and $V_s$ are uniformly bounded on $D_{2K_0s}$.
Lemma~\ref{lem:real-complex} makes their difference small on
$D_{K_0s}$, with a fixed positive power of $h_M(C_*s)$.
Lemma~\ref{lem:weight-loss} absorbs that power, proving
\eqref{pos:central-bounds}.
\end{proof}

\subsection{Cauchy division and quantitative remainders}

The integral in the next lemma is the classical Cauchy construction
in holomorphic Weierstrass division; see
\cite[Theorem~6.2.5]{LeblSCV}. We state its quantitative form for a
contour that may enclose only some roots. The bounds follow from
contour separation and the residue formula; no division hypothesis
on $H$ is required.

\begin{lemma}[Cauchy division estimates]\label{app:cauchy-division}
Let $\mathcal Z\subset\C$ be a disc, let $W(z,w)$ be monic of
degree $m\ge1$ in $w$, with coefficients holomorphic on
$\mathcal Z$, and let $\Gamma$ be a positively oriented,
piecewise $C^1$ Jordan curve. Write $\mathcal V$ for its interior
and let $\Omega\subset\mathcal V$ be open. Suppose that $H$ is
holomorphic on a neighbourhood of
$\mathcal Z\times\overline{\mathcal V}$ and that
$W(z,\zeta)\ne0$ for $z\in\mathcal Z$, $\zeta\in\Gamma$.
Then
\begin{equation}\label{app:cauchy-integral}
 G(z,w)=\frac{1}{2\pi i}\int_\Gamma
  \frac{H(z,\zeta)}{W(z,\zeta)(\zeta-w)}\,d\zeta
\end{equation}
is holomorphic on $\mathcal Z\times\mathcal V$.

Suppose, for $0<\sigma<1$, that the length of $\Gamma$ is at most
$C_0\sigma$, that $\operatorname{dist}(\Omega,\Gamma)\ge c_0\sigma$,
and that every root of $W(z,\cdot)$ has distance at least
$c_0\sigma$ from $\Gamma$, uniformly in $z$. Then
\begin{equation}\label{app:cauchy-norm}
 \|G\|_{\mathcal Z\times\Omega}
 \le C\sigma^{-m}\|H\|_{\mathcal Z\times\Gamma}.
\end{equation}

If the roots $\lambda_1,\ldots,\lambda_m$ are distinct and
holomorphic on $\mathcal Z$, let $I$ be the set of roots enclosed
by $\Gamma$. Then
\begin{equation}\label{pos:residue-remainder}
 H(z,w)-W(z,w)G(z,w)
 =\sum_{i\in I}
 \frac{H(z,\lambda_i(z))}{\partial_wW(z,\lambda_i(z))}
 \frac{W(z,w)}{w-\lambda_i(z)}.
\end{equation}
If $|\lambda_i(z)-\lambda_j(z)|\ge c_1\sigma^{A_0}$ for $i\ne j$,
if the roots and $\Omega$ stay in a fixed bounded neighbourhood,
and if $|H(z,\lambda_i(z))|\le\varepsilon$ for $i\in I$, then
\begin{equation}\label{app:cauchy-remainder}
 \|H-WG\|_{\mathcal Z\times\Omega}
 \le C\sigma^{-A_0(m-1)}\varepsilon.
\end{equation}

Finally, without any simplicity assumption on the roots, suppose
$H=WQ+E$, with $Q,E$ holomorphic on a neighbourhood of
$\mathcal Z\times\overline{\mathcal V}$. Under the same contour
separation assumptions,
\begin{equation}\label{app:cauchy-collision}
 \|G-Q\|_{\mathcal Z\times\Omega}
 \le C\sigma^{-m}\|E\|_{\mathcal Z\times\Gamma}.
\end{equation}
If $W$ is uniformly bounded on $\mathcal Z\times\Omega$, this gives
\[
 \|H-WG\|_{\mathcal Z\times\Omega}
 \le (1+C\sigma^{-m})
       \|E\|_{\mathcal Z\times\overline{\mathcal V}}.
\]
Constants depend only on the fixed separation, length, degree, and
boundedness constants.
\end{lemma}

\begin{proof}
The integrand in \eqref{app:cauchy-integral} is holomorphic in
$(z,w)$ off the fixed contour. Differentiation under the integral
therefore gives holomorphicity. Factoring the monic polynomial
shows that $|W(z,\zeta)|\ge(c_0\sigma)^m$ on $\Gamma$.
Together with its length and the distance of $\Omega$ from
$\Gamma$, this proves \eqref{app:cauchy-norm}.

For fixed $z$ and $w$ away from the roots, the residue theorem gives
the contribution $H(z,w)/W(z,w)$ from the pole at $\zeta=w$, and
the contribution
\[
 \frac{H(z,\lambda_i(z))}
 {\partial_wW(z,\lambda_i(z))(\lambda_i(z)-w)}
\]
from each enclosed root. Multiplying by $W(z,w)$ gives
\eqref{pos:residue-remainder}. Each
$W(z,w)/(w-\lambda_i(z))$ is a polynomial in $w$, so the formula
extends holomorphically across the roots. The set $I$ is independent
of $z$, since a root cannot cross $\Gamma$. Furthermore,
\[
 |\partial_wW(z,\lambda_i(z))|
 =\prod_{j\ne i}|\lambda_i(z)-\lambda_j(z)|
 \ge c_1^{m-1}\sigma^{A_0(m-1)}.
\]
The numerator polynomials $W(z,w)/(w-\lambda_i(z))$ are uniformly
bounded under the stated boundedness assumptions. Summing over at
most $m$ enclosed roots proves \eqref{app:cauchy-remainder}.

If $H=WQ+E$, the contribution of $WQ$ to the integral is $Q$ by
Cauchy's formula. Applying \eqref{app:cauchy-norm} to $E$ proves
\eqref{app:cauchy-collision}. The last assertion follows from
$H-WG=E-W(G-Q)$.
\end{proof}

\subsection{Removal of a small polynomial factor}

\paragraph{Finite differences.}\label{app:small-product-proof}
Lemma~\ref{pos:small-product} is an elementary consequence of the
$m$th finite-difference identity for a monic polynomial of degree
$m$. A first derivative bound transfers a value estimate from a
nearby point to the point where $W$ may vanish.

\begin{proof}[Proof of Lemma~\ref{pos:small-product}]
Fix $(x,b)$ in the smaller square and $0<\tau\le r/2$.
The $m$th finite difference of $t\mapsto W(x,b+t)$ at the
$m+1$ equally spaced points of $[-\tau,\tau]$ equals
$m!(2\tau/m)^m$. The sum of the absolute values of its coefficients
is $2^m$, so $|W(x,b+t)|\ge m!\tau^m/m^m$ at one of these points.
The product and derivative bounds give
\[
 |v(x,b)|\le C_m\varepsilon\tau^{-m}+L_0\tau.
\]
Take $\tau=\varepsilon^{1/(m+1)}$.
\end{proof}
\section{Proofs for the sharpness construction}
\label{app:sharpness}

We prove the statements of Sections~\ref{sec:staircase}--\ref{sec:counterexample}
in dependency order, with the factorial kept outside the weight.

\subsection{The staircase weight}\label{app:stair-details}
For the standard Denjoy--Carleman criteria, see
\cite{Mandelbrojt1952,Rainer2022Survey}.

\begin{proof}[Proof of Lemma~\ref{stair:separation}]
Put $x_j=E_{6j}(2)$.  Then $r_j=x_j+O(1)$,
$s_j=E_4(x_j)+O(1)$, and $r_{j+1}=E_6(x_j)+O(1)$.  Hence
\[
 \frac{\log((r_{j+1}+1)/(s_j+1))}{Q(s_j)}
 \sim\frac{E_5(x_j)}{E_3(x_j)^\alpha E_2(x_j)^{2\alpha}}
 \longrightarrow\infty.
\]
For fixed $C>0$, eventually $Cr_j\leq e^{x_j}/2$, so
\[
 \log\frac{E_3(Cr_j)}{E_4(x_j)}
 \leq E_2(e^{x_j}/2)-E_2(e^{x_j})\longrightarrow-\infty.
\]
This proves \eqref{stair:degree-room}.  In particular,
$s_j\geq E_3(r_j)$ eventually, whence
\[
 \frac{Q(r_j)^{1/\alpha}}{(\log L(s_j))^2}
 \leq Q(r_j)^{1/\alpha}e^{-2r_j}\longrightarrow0.
\]
Finally, $\log r_j\geq E_{6j-1}(2)\geq j$.
\end{proof}

\begin{proof}[Proof of Proposition~\ref{stair:regularity}]
In the coordinate $u=\log Q(k)$, the staircase has slopes zero and
one, starts at height zero, and stays below $u$.  Thus
$1\leq\mu_k\leq Q(k)$, the quotients are nondecreasing, and
$\mu_k/Q(k)$ is nonincreasing.  Formula~\eqref{stair:cj} gives
$c_j\geq Q(r_j)^{-1}$, while $\mu_k\leq Q(k)$ on the $j$th
growth interval gives $c_j\leq1$.  This proves \eqref{stair:bounds}
and the monotonicity and logarithmic convexity of $M$.

Integrating
\[
 \frac{d\log Q(t)}{d\log t}
 =\frac{\alpha t}{(t+e^e)L(t)}
       \left(1+\frac2{\log L(t)}\right)\leq3\alpha
\]
gives
\begin{equation}\label{stair:power-bound}
 \frac{Q(t)}{Q(u)}\leq(t/u)^{3\alpha}\qquad(t\geq u\geq1).
\end{equation}
Consequently, for integers $p,\ell\geq1$,
\[
 \frac{M_{p+\ell}}{M_pM_\ell}
 =\prod_{i=1}^\ell\frac{\mu_{p+i}}{\mu_i}
 \leq\prod_{i=1}^\ell\frac{Q(p+i)}{Q(i)}
 \leq\binom{p+\ell}{\ell}^{3\alpha}\leq2^{3\alpha(p+\ell)}.
\]
The cases $p=0$ or $\ell=0$ are immediate.  This proves moderate
growth; taking $\ell=1$ and using $M_1=1$ also gives differentiation
stability.

On each plateau $\mu_k=c_jQ(s_j)\leq Q(s_j)$, so
\[
 \sum_{k=s_j+1}^{r_{j+1}}\frac1{k\mu_k}
 \geq\frac1{Q(s_j)}\log\frac{r_{j+1}+1}{s_j+1}
 \longrightarrow\infty.
\]
The disjoint contributions prove the quasianalyticity sum diverges.
Also $\mu_{s_j}\geq Q(s_j)/Q(r_j)\to\infty$, and monotonicity
gives $\mu_k\to\infty$.  Since $k^{-1}\log M_k$ is the arithmetic
mean of $\log\mu_1,\ldots,\log\mu_k$, it tends to infinity.
The weight comparison criterion gives strict inclusion of the analytic
class \cite[\S1.4]{Thilliez2008}; the remaining closure properties
follow from the usual calculus
\cite[\S\S2--4]{BierstoneMilman2004}, \cite{Thilliez2008}.
\end{proof}

\subsubsection{Quotient counting and the growth threshold}
\label{app:stair-counting}

\begin{proof}[Verification of \eqref{stair:omega-count}]
Let $N(T)=\#\{k\geq1:\mu_k\leq T\}$, which is finite.  The
successive increments of $k\log T-\log M_k$ are
$\log(T/\mu_k)$: they are nonnegative through $N(T)$ and negative
afterwards.  The maximum is therefore the sum of those first
$N(T)$ increments, including any zero increments.  This proves the
formula, also when $N(T)=0$.
\end{proof}

\begin{proof}[Proof of Proposition~\ref{stair:threshold}]
\emph{Lower scale growth.}
Since $L(k^2)\leq2L(k)$ and $\log L(k)\geq1$,
\[
 Q(k^2)\leq DQ(k),\qquad D=2^\alpha(1+\log2)^{2\alpha}.
\]
The monotonicity of $\mu_k/Q(k)$ implies
\begin{equation}\label{stair:square-count}
 \mu_{k^2}\leq D\mu_k.
\end{equation}
Set $R=2D$ and $N=N(T/R)$.  Then $\mu_{N^2}\leq T/2$, so
\eqref{stair:omega-count} gives
\[
 \omega_M(T)\geq N^2\log2,\qquad
 \omega_M(T/R)\leq N\log T.
\]
Taking logarithms yields
\begin{equation}\label{stair:lambda-count}
 \Lambda_M(T)-\Lambda_M(T/R)
 \geq\log N(T/R)-\log\log T+\log\log2.
\end{equation}
For a sufficiently small fixed $d>0$, put
$u=dT^{1/\alpha}/(\log T)^2$.  Then
\[
 u^\alpha(\log u)^{2\alpha}\sim d^\alpha\alpha^{-2\alpha}T<T/R.
\]
Every integer $k\leq e^u-e^e$ satisfies $\mu_k\leq Q(k)\leq T/R$.
Thus
\begin{equation}\label{stair:count-lower}
 \log N(T/R)\geq c_1\frac{T^{1/\alpha}}{(\log T)^2}
\end{equation}
for large $T$.  Combining this with \eqref{stair:lambda-count} proves
\eqref{stair:log-threshold}.

\emph{Failure at the endpoint.}
Set $T_j=\mu_{s_j}=c_jQ(s_j)\to\infty$.  The following plateau
contributes zero to \eqref{stair:omega-count}, giving
\[
 \omega_M(T_j)\leq s_j\log T_j,\qquad
 \Lambda_M(T_j)\leq C\log s_j,
\]
where the second inequality uses $T_j\leq Q(s_j)$.  Therefore
\begin{equation}\label{stair:endpoint-failure}
 \frac{\Lambda_M(T_j)}{T_j^{\beta_0}}
 \leq\frac{C c_j^{-1/\alpha}}{(\log L(s_j))^2}
 \leq\frac{C Q(r_j)^{1/\alpha}}{(\log L(s_j))^2}
 \longrightarrow0.
\end{equation}
For each fixed $R'>1$, eventually $\Lambda_M(T_j/R')\geq0$, so the
same limit holds for
$[\Lambda_M(T_j)-\Lambda_M(T_j/R')]/T_j^{\beta_0}$.
No choice of $R'$ and positive constant gives $(G_{\beta_0})$.
Finally, \eqref{stair:log-threshold} gives $(G_\beta)$ for every
$\beta<\beta_0$, whereas any $(G_\beta)$ with $\beta>\beta_0$
would imply $(G_{\beta_0})$.  This proves \eqref{stair:spectrum}.
\end{proof}

\subsection{The auxiliary associated function}\label{app:products-cutoffs}

\begin{proof}[Proof of Lemma~\ref{neg:product-asymptotics}]
For $a>0$, $b\geq0$, consider
$S_k=\prod_{i=1}^kL(i)^a(\log L(i))^b$.  Let $u=u(T)>e$ solve
$u^a(\log u)^b=T$.  Its quotient-counting function is
$N_S(T)=\lfloor e^u-e^e\rfloor$.  The counting formula gives
\[
 \omega_S(T)\leq e^u\log T.
\]
For $i\leq\lfloor e^{u-1}-e^e\rfloor$, each summand is at least
$a\log(u/(u-1))\geq a/u$.  There are at least $ce^u$ such indices,
so
\[
 ce^u/u\leq\omega_S(T)\leq e^u\log T,\qquad
 \log\omega_S(T)=u+O(\log u).
\]
Elementary inversion gives
\[
 u\sim a^{b/a}T^{1/a}(\log T)^{-b/a}.
\]
Indeed the expression on the right, denoted $u_0$, satisfies
$u_0^a(\log u_0)^b\sim T$; substituting $(1\pm\varepsilon)u_0$
and using monotonicity proves the assertion.  Taking
$a=\alpha-1$, $b=2\alpha-3/2$ proves \eqref{neg:D-asymptotics}.
Moreover, $\log\log\omega_{\mathsf D}(T)\sim
(\alpha-1)^{-1}\log T$, and hence
\[
 Q(\omega_{\mathsf D}(T))
 \asymp T^{\alpha/(\alpha-1)}(\log T)^{-\alpha\nu+2\alpha}
 =T^{q_*}(\log T)^{-q_*/2}.
\]
This is \eqref{neg:composed-asymptotic}.
\end{proof}

\subsection{The localized smooth model}\label{app:localized-blocks}

\begin{proof}[Proof of Lemma~\ref{neg:localized-estimates}]
\emph{The cutoff.}
Put $b_k=\mathsf B_k/\mathsf B_{k-1}$ and
$S=\sum_{k\geq1}(kb_k)^{-1}<\infty$.  Set
$d_k=(16Skb_k)^{-1}$ and
$H_d=d^{-1}\mathbf1_{[-d/2,d/2]}$.  The infinite-convolution
construction \cite[Lemma~3.3]{Rainer2022Survey} gives the smooth
probability density $\Phi=*_{k\geq1}H_{d_k}$, supported in
$[-1/32,1/32]$, with
\[
 \|\Phi^{(r)}\|_\infty
 \leq\frac{2^r}{d_1\cdots d_r d_{r+1}}
 =2^r(16S)^{r+1}(r+1)!\mathsf B_{r+1}
 \leq C^{r+1}r!\mathsf B_r.
\]
The last inequality uses $(r+1)b_{r+1}\leq C_1^{r+1}$.
Then $\chi=\mathbf1_{[-1/8,1/8]}*\Phi$ has support in
$[-5/32,5/32]$, takes values in $[0,1]$, equals one near zero,
and satisfies \eqref{neg:cutoff-bounds}.

\emph{The rational factor.}
In the variables $\xi=(x-a)/a$, $\eta=y/a^q$, the reciprocal is
$a^{-2n}/F(\xi,\eta)$, where
$F(\xi,\eta)=(1+\xi)^{2n}+\eta^{2m}$.
For real $\xi\in[-1/3,1/3]$, $t\in\R$, and complex
$|z|,|w|\leq\rho\leq1$, the polynomial formula gives
\[
 F(\xi,t)\geq c(1+|t|^{2m}),\qquad
 |F(\xi+z,t+w)-F(\xi,t)|\leq C\rho(1+|t|^{2m}).
\]
A sufficiently small fixed $\rho$ therefore gives a zero-free
polydisc and a uniform bound for $1/F$ at every real centre.
Cauchy's estimates, followed by Leibniz's formula with $\chi(\xi)$,
give $C^{k+1}k!\mathsf B_k$ in the scaled variables.  Each derivative
in the original variables costs at most $a^{-q}$, proving
\eqref{neg:quotient-bounds}.  These are the same type of rational-block
estimates as in \cite[\S3.2, Lemma~3.3]{BelottoBierstoneKiro2022}.
The support stays away from $x=0$, so extension by zero is smooth.
Finally, $Pv_a=\chi((x-a)/a)$ proves \eqref{neg:numerator-bounds},
and the geometric series at $x=a$ gives
\eqref{neg:exact-derivative}.
\end{proof}

\subsection{Temporary weights and the weighted smooth blocks}
\label{app:temporary-weights}\label{app:weighted-blocks}

\begin{proof}[Proof of Lemma~\ref{neg:temporary-properties}]
The ratio $\mu_i/Q(i)$ is nonincreasing and equals $c_j$ for
$r_j\leq i\leq s_j$.  Hence the new quotients equal $\mu_i$
through $s_j$ and $c_jQ(i)$ afterwards.  They are nondecreasing and
at least one, and their products agree with $M$ through $s_j$.
Each quotient is also at least $c_jQ(i)$, giving
\eqref{neg:temporary-bounds}.  The tail sum
$\sum_i(ic_jQ(i))^{-1}$ converges because $\alpha>1$; the
temporary class is therefore nonquasianalytic.
\end{proof}

\begin{proof}[Proof of Lemma~\ref{neg:large-derivative}]
For $k=|\gamma|$, the numerator and temporary-weight bounds give
\[
 \frac{\|\partial^\gamma(Pu_j)\|_K}{k!\widetilde M^{(j)}_k}
 \leq w_jC_\chi\frac{\tau_j^k}{\mathsf D_k}
 \leq w_jC_\chi e^{\omega_{\mathsf D}(\tau_j)}=\eta_j.
\]
This proves \eqref{neg:small-numerator} at every order.
Write $r=r_j$ and $c=c_j$.  Since
$1\leq c^{-1}\leq Q(r)$, we have $\log\tau_j\sim r$.
The derivative bound gives
$\log\|v_{a_j}\|_{C^r(K)}\leq Cr^2$, whereas
\eqref{neg:D-asymptotics} gives
\begin{equation}\label{neg:amplitude-decay}
 \omega_{\mathsf D}(\tau_j)/r^2\longrightarrow\infty.
\end{equation}
The definition of $w_j$ therefore proves
\eqref{neg:small-low-orders}.

The prescribed order satisfies
\begin{equation}\label{neg:order-size}
 k_j=\log(1/w_j)+O(1)
     =\omega_{\mathsf D}(\tau_j)+O(j+1)
     \sim\omega_{\mathsf D}(\tau_j).
\end{equation}
In particular $k_j>r_j$ for large $j$.  The upper estimate in
\eqref{neg:D-asymptotics} also gives
$\log k_j\leq e^{Cr_j}$ for a fixed $C$, and thus
$k_j\leq E_2(Cr_j)<s_j$ by \eqref{stair:degree-room}.

By construction $w_j^{1/k_j}\geq e^{-1}$.  Since
$M_k\leq Q(k)^k$, the exact derivative gives
\begin{equation}\label{neg:amplitude-compensation}
 \left(\frac{|\partial_y^{k_j}u_j(z_j)|}{k_j!M_{k_j}}\right)^{1/k_j}
 \geq\frac{e^{-1}a_j^{-q}}{Q(k_j)}.
\end{equation}
Here the additional factor $a_j^{-2n/k_j}\geq1$ was discarded.
By \eqref{neg:order-size} and \eqref{neg:composed-asymptotic},
\[
 Q(k_j)\asymp Q(\omega_{\mathsf D}(\tau_j))
 \leq Cc^{-q_*}e^{q_*r}r^{-q_*/2}.
\]
Substitution in \eqref{neg:amplitude-compensation}, followed by
$c\geq Q(r)^{-1}$, proves \eqref{neg:endpoint-gap}.
Finally, because $q\geq q_*$ and $j\leq\log r_j$,
\begin{equation}\label{neg:decisive-comparison}
 \frac1j e^{(q-q_*)r_j}
 \left(\frac{\sqrt{r_j}}{Q(r_j)}\right)^{q_*}
 \geq\frac{r_j^{q_*/2}}{(\log r_j)Q(r_j)^{q_*}}
 \longrightarrow\infty.
\end{equation}
Thus the right side of \eqref{neg:endpoint-gap} is eventually at
least $8j$, proving \eqref{neg:large-derivative-bound}, including
the endpoint.
\end{proof}

\subsection{Quantitative Gaussian--Taylor approximation}
\label{app:negative-approx}

Prescribe
\begin{equation}\label{neg:approximation-prescription}
 \kappa=2+\frac{2q}{\alpha-1},\qquad
 \epsilon_j=\exp\bigl(-E_2(\kappa r_j)\bigr),\qquad
 d_j=\left\lceil\epsilon_j^{-3}\right\rceil.
\end{equation}
Let $G_\epsilon(t)=(2\pi\epsilon^2)^{-1}
\exp(-|t|^2/(2\epsilon^2))$ on $\R^2$.
The polynomial $p_j$ is the Taylor polynomial of
$U_j=u_j*G_{\epsilon_j}$ at the origin, truncated at degree $d_j$
separately in each variable.

\begin{proof}[Proof of Lemma~\ref{neg:gaussian-approximation}]
\emph{Smoothing in the weighted norm.}
Write $r=r_j$, $c=c_j$, $S=\widetilde M^{(j)}$.
Since $w_j\leq1$, the global derivative bounds imply
\[
 \|\partial^\gamma u_j\|_\infty
 \leq C_v^{k+1}e^{2nr+qrk}k!\mathsf B_k,
 \qquad k=|\gamma|.
\]
Also $(k+1)\mathsf B_{k+1}\leq C_B^{k+1}\mathsf B_k$.
The mean value theorem and the first absolute moments of the
Gaussian, followed by division by $k!S_k\geq k!c^k\mathsf A_k$,
therefore give the single all-order estimate
\begin{equation}\label{neg:weighted-smoothing}
 \|u_j*G_\epsilon-u_j\|_{S,K}
 \leq C\epsilon\exp\bigl((2n+q)r+
              \omega_{\mathsf D}(Ce^{qr}/c)\bigr).
\end{equation}
For example, before division the derivative of order $k$ is bounded
by $C\epsilon C_v^{k+2}e^{(2n+q)r+qrk}(k+1)!\mathsf B_{k+1}$;
taking the supremum produces \eqref{neg:weighted-smoothing}.
By \eqref{neg:D-asymptotics} and $c^{-1}\leq Q(r)$,
\[
 \log\omega_{\mathsf D}(Ce^{qr}/c)
 \leq\exp\bigl((q/(\alpha-1)+1)r\bigr)
\]
for large $r$.  Since $\kappa>q/(\alpha-1)+1$, we have
$\omega_{\mathsf D}(Ce^{qr}/c)=o(E_2(\kappa r))$.
With $\epsilon=\epsilon_j$, we obtain
\begin{equation}\label{neg:appendix-smoothing}
 \|U_j-u_j\|_{S,K}
 \leq\exp\bigl(-\tfrac12E_2(\kappa r)\bigr)
 \leq\tfrac12e^{-r^2}.
\end{equation}

\emph{Taylor truncation.}
The Gaussian convolution of the bounded function $u_j$ extends
entirely by its locally uniformly convergent integral and satisfies
\[
 |U_j(z)|\leq\|u_j\|_\infty
 \exp\left(\frac{|\operatorname{Im}z|^2}{2\epsilon_j^2}\right).
\]
Hence $\log\|U_j\|_\infty\leq C\epsilon_j^{-2}$ on
$|z_1|,|z_2|\leq8$.  Cauchy's coefficient estimates and two
geometric series give, on $|z_1|,|z_2|\leq4$,
\[
 \|U_j-p_j\|_\infty
 \leq8\exp(C\epsilon_j^{-2})2^{-d_j}
 \leq\tfrac12e^{-r^2},
\]
because $d_j\geq\epsilon_j^{-3}$.  Each point of $K$ is the centre
of a complex polydisc of radius $2$ contained in this latter
polydisc.  Another Cauchy estimate gives, for every $\gamma$,
\[
 \|\partial^\gamma(U_j-p_j)\|_K
 \leq\tfrac12e^{-r^2}\gamma!2^{-|\gamma|}
 \leq\tfrac12e^{-r^2}|\gamma|!S_{|\gamma|}.
\]
Together with \eqref{neg:appendix-smoothing}, this proves the norm
bound.  Real-valued $u_j$ gives real Taylor coefficients.

\emph{Degree.}
By construction $\deg p_j\leq2d_j$ and
$\log d_j\leq4E_2(\kappa r_j)\leq E_2(Cr_j)$ for any fixed
$C>\kappa$ and all sufficiently large $j$.  Thus
$d_j\leq E_3(Cr_j)$, and \eqref{stair:degree-room} gives
$2d_j+2n<s_j$.
\end{proof}

\subsection{Transfer to the polynomial blocks}\label{app:subsidiary-blocks}

\begin{proof}[Verification of \eqref{neg:polynomial-multiplication}]
For $k=|\gamma|$, $l=|\delta|$,
$\binom\gamma\delta(k-l)!/k!\leq1/\delta!$ and $S_{k-l}\leq S_k$.
Leibniz's formula therefore gives
\[
 \frac{\|\partial^\gamma(Fu)\|_K}{k!S_k}
 \leq\sum_{\delta\leq\gamma}
       \frac{\|\partial^\delta F\|_K}{\delta!}\|u\|_{S,K}.
\]
Taking the supremum proves the estimate.
\end{proof}

\begin{proof}[Proof of Proposition~\ref{neg:polynomial-blocks}]
The order and degree bounds are already established.  Equality of
the weights through $r_j$, together with $M_{r_j}\leq Q(r_j)^{r_j}$,
gives
\[
 \|p_j-u_j\|_{C^{r_j}(K)}
 \leq\delta_j r_j!M_{r_j}
 \leq\exp\bigl(-r_j^2+r_j\log r_j+r_j\log Q(r_j)\bigr)
 \leq\eta_j
\]
for large $j$.  Thus $\|p_j\|_{C^{r_j}(K)}\leq2\eta_j$.
Multiplication and \eqref{neg:finite-degree-transfer} give
\[
 \|Pp_j\|_{M,K}
 =\|Pp_j\|_{\widetilde M^{(j)},K}
 \leq\eta_j+C_P\delta_j\leq2\eta_j.
\]
These also prove \eqref{neg:polynomial-transfer}; their sum is
at most $4\eta_j\leq2^{-j}$.
At $k_j<s_j$, the weights still agree, so
\[
 \frac{|\partial_y^{k_j}p_j(z_j)|}{k_j!M_{k_j}}
 \geq(8j)^{k_j}-\delta_j\geq(4j)^{k_j}.
\]
This is the last assertion of the proposition.
\end{proof}

\subsection{Summation and the counterexample}\label{app:negative-summation}

\begin{proof}[Proof of Theorem~\ref{thm:sharpness}]
Choose $J$ so that Proposition~\ref{neg:polynomial-blocks} holds for
$j\geq J$, and define $g,f$ by \eqref{neg:final-functions}.
For each fixed $r$, eventually $r\leq r_j$, so
$\|p_j\|_{C^r(K)}\leq2^{-j}$.  Every derivative series converges
uniformly, giving smooth functions on $\operatorname{int}K$.
The numerator bounds imply $f\in\CM(\operatorname{int}K)$ and,
on $K_0=[-1/2,1/2]^2$,
\begin{equation}\label{neg:fixed-norm-convergence}
 \left\|f-\sum_{j=J}^NPp_j\right\|_{M,K_0}
 \leq\sum_{j>N}2^{-j}\longrightarrow0.
\end{equation}
Indeed, the same bound holds at every normalized derivative order
before taking the supremum.

At order $k_j$, all earlier polynomials vanish, since
$\deg p_i<s_i<r_j<k_j$ for $i<j$.  For $i>j$ we have
$k_j<s_j<r_i$, so the low-order bound alone gives
\begin{equation}\label{neg:block-later-errors}
 \frac{|\partial_y^{k_j}p_i(z_j)|}{k_j!M_{k_j}}
 \leq\frac{\|p_i\|_{C^{r_i}(K)}}{k_j!M_{k_j}}\leq2^{-i}.
\end{equation}
Therefore
\begin{equation}\label{neg:noncancellation}
 \frac{|\partial_y^{k_j}g(z_j)|}{k_j!M_{k_j}}
 \geq(4j)^{k_j}-\sum_{i>j}2^{-i}\geq(2j)^{k_j}.
\end{equation}
As $z_j\to0$ and $k_j\to\infty$, this contradicts every bound
$AB^{k_j}$ on a compact neighbourhood of the origin.
Thus $g$ belongs to $\CM$ on no such neighbourhood.
Finally, $P>0$ off the origin, so any smooth quotient of $f$ by $P$
equals $g$ there and hence everywhere by continuity.  Thus
$f\notin P\CM$, while \eqref{neg:fixed-norm-convergence} supplies
the asserted fixed-norm approximation by polynomial multiples.
\end{proof}

\section{Proof of the Division by quasianalytic germs and Noetherianity}

\subsection{Proof of Lemma \ref{lem:reduced-normal-form}}\label{app:reduced-normal-form}

We first record some elementary facts about powers of the maximal ideal. Since $\CM$ is stable under division by a coordinate,
for every $k\geq1$ we have
\begin{equation}\label{eq:algebra-jets}
\widehat f\in\mathfrak n^k
\quad\Longleftrightarrow\quad
f\in\mathfrak m^k,
\qquad
A/\mathfrak m^k\simeq B/\mathfrak n^k.
\end{equation}
In particular, $\mathfrak m^k$ is generated by the monomials of degree $k$, and quasianalyticity gives
\[
\bigcap_{k\geq1}\mathfrak m^k=0.
\]
We are ready to prove:

\begin{proof}[Proof of Lemma \ref{lem:reduced-normal-form}]
Since $\widehat \varphi$ is formally reduced, its partial derivatives
have no common irreducible factor. As $B$ is a
two-dimensional regular local ring, it follows that
\begin{equation}
 \mathfrak n^k\subset(\widehat \varphi_x,\widehat \varphi_y)B
 \quad\hbox{for some }k\geq1.
 \label{eq:formal-jacobian-colength}
\end{equation}
Now, let $J=(\varphi_x,\varphi_y)A$; we need to show the analogous statement for $J$. Indeed, by \eqref{eq:formal-jacobian-colength} and
\eqref{eq:algebra-jets}, we have $\mathfrak m^k\subset J+\mathfrak m^{k+1}$. Now, since $\mathfrak m^k$ is finitely generated,
the module
\[
L=(\mathfrak m^k+J)/J
\]
is finitely generated and satisfies $L=\mathfrak mL$. Nakayama's lemma (which does not require the ambient ring to be Noetherian, see~\cite[Tag~00DV]{StacksProject}) gives
\begin{equation}
 L=0 \implies \mathfrak m^k\subset J.
 \label{eq:actual-jacobian-colength}
\end{equation}
Next, choose $N\geq k+2$, let $P$ be the Taylor polynomial of $\varphi$ of degree
$N$, and put $R=\varphi-P$.  Equations~\eqref{eq:algebra-jets}
and~\eqref{eq:actual-jacobian-colength} give
\[
 R=b^T\nabla \varphi,\qquad \nabla R=C\nabla \varphi,
 \qquad b\in(\mathfrak m^3)^2,
 \quad C\in\operatorname{Mat}_2(\mathfrak m^2),
\]
where gradients are columns.  Set $\varphi_t=\varphi-tR$ and
\[
 V_t=(I-tC)^{-T}b,\qquad 0\leq t\leq1.
\]
On a sufficiently small fixed neighbourhood, $I-tC$ is invertible for
all $t$.  The vector field $V_t$ is of class $\CM$ and vanishes to order at least three at the origin.  Its flow $\Psi_t$, starting at
the identity, exists for $0\leq t\leq1$ after a further restriction of the neighbourhood and preserves $\CM$ regularity by
\cite[Theorem~1 and Remark~4]{RainerSchindl2016}. 

We now show that $\varphi=P\circ\Psi_1$; note that it is enough to show that the derivative of $\varphi_t(\Psi_t)$ is $0$. Indeed, since
$\nabla \varphi_t=(I-tC)\nabla \varphi$, we have
\[
 \frac{d}{dt}\varphi_t(\Psi_t)
 =-R(\Psi_t)+\nabla \varphi_t(\Psi_t)\cdot V_t(\Psi_t)=0.
\]
Consequently $\varphi=P\circ\Psi_1$ as needed. Finally, the induced formal coordinate change preserves squarefreeness, since it is a formal diffeomorphism, so $P$ is reduced.
\end{proof}
\section[Failure of division in a Gevrey quasianalytic class]
{Failure of division in a Gevrey quasianalytic class---construction by ChatGPT}
\label{app:sectorial-gevrey}

We give a counterexample to exact division in a sectorial quasianalytic
class. Here \emph{Gevrey order one} refers to the summability convention:
Taylor coefficients have factorial growth, and derivatives have squared
factorial bounds. The example concerns germs on the \emph{positive
quadrant}; its defining integral is not defined for negative $x$.

Put
\[
 S_\theta(r)=\{z\in\mathbb C:0<|z|<r,\ |\arg z|<\theta\}.
\]
For $d=1,2$, let $\mathcal S_1^{(d)}$ consist of germs that are
$C^\infty$ up to the boundary of $[0,r)^d$ and admit a holomorphic
extension $U$ to $S_\theta(r)^d$ for some $\theta>\pi/2$, such that,
on every smaller product $S_{\theta'}(r')^d$, there are $C,A>0$ with
\begin{equation}\label{eq:D-sectorial-bound}
 \sup_{S_{\theta'}(r')^d}|\partial^\alpha U|
 \le C A^{|\alpha|}(|\alpha|!)^2
 \qquad(\alpha\in\mathbb N^d).
\end{equation}
Thus each sector has opening $2\theta>\pi$. In one variable these
are the usual Gevrey-one sums in direction zero; see
\cite{Balser1994}. Watson's uniqueness theorem
\cite[\S8]{Watson1912} makes this class quasianalytic: for a flat
germ, apply the theorem to $t\mapsto U(ta)$ for each
$a\in(0,\infty)^d$. We use precisely this product-sector class.

\begin{proposition}\label{prop:D-cusp-summability}
The function
\begin{equation}\label{eq:D-numerator}
 \boxed{\displaystyle
 f(x,y)=\int_0^\infty e^{-t}
 \left(\frac1{1+y^2t^2}-\frac1{1+x^3t^2}\right)\,dt}
 \qquad(x,y\ge0)
\end{equation}
belongs to $\mathcal S_1^{(2)}$. It is divisible by $P=y^2-x^3$
with a quotient $g$ that is $C^\infty$ up to the boundary of the
positive quadrant, but $g\notin\mathcal S_1^{(2)}$.
\end{proposition}

\begin{proof}
Let
\[
 H(z)=\int_0^\infty\frac{e^{-t}}{1+z^2t^2}\,dt
 \qquad(z>0).
\]
Its formal expansion is
$\widehat H(z)=\sum_{n\ge0}(-1)^n(2n)!z^{2n}$, whose
order-one Borel transform is $(1+\xi^2)^{-1}$. Rotating the
Laplace ray gives its holomorphic continuation to $|\arg z|<\pi$:
\begin{equation}\label{eq:D-laplace}
 H(z)=\frac1z\int_0^{e^{i\phi}\infty}
       \frac{e^{-\xi/z}}{1+\xi^2}\,d\xi,
 \qquad |\phi|<\frac\pi2,\quad
        |\arg z-\phi|<\frac\pi2.
\end{equation}
The geometric-series remainder gives Gevrey-one asymptotics on
every smaller sector. Hence $H(y)$ has order one, whereas
$H(x^{3/2})$ has order $2/3$, on sectors
$|\arg x|<2\pi/3$; the power uses the branch positive for $x>0$.
Consequently $f=H(y)-H(x^{3/2})$ satisfies
\eqref{eq:D-sectorial-bound} on smaller products whenever
$\pi/2<\theta<2\pi/3$.

Direct subtraction in \eqref{eq:D-numerator} gives $f=Pg$, where
\begin{equation}\label{eq:D-quotient}
 g(x,y)=-\int_0^\infty
 \frac{e^{-t}t^2}{(1+y^2t^2)(1+x^3t^2)}\,dt.
\end{equation}
This formula also defines $g$ on the real cusp. All derivatives
may be passed under the integral, including at the boundary.
Indeed, Cauchy's estimates for $(1+s^3)^{-1}$ and $(1+s^2)^{-1}$
on a fixed complex neighbourhood of $[0,\infty)$ give
\begin{equation}\label{eq:D-real-bound}
 |\partial_x^p\partial_y^q g(x,y)|
 \le C A^{p+q}p!q!\,
       \Gamma\!\left(3+\frac{2p}{3}+q\right)
 \le C B^{p+q}\bigl((p+q)!\bigr)^2.
\end{equation}

Nevertheless, $g$ has no holomorphic extension to any
$S_\theta(r)^2$ with $\theta>\pi/2$. To see this, rotating the
contours in \eqref{eq:D-laplace} past the pole $\xi=i$ yields
\begin{equation}\label{eq:D-stokes-jump}
 H(w)-H(-w)=\frac\pi w e^{-i/w},
 \qquad 0<\arg w<\pi,
\end{equation}
where $\arg(-w)=\arg w-\pi$. For small $\rho>0$, put
\[
 x=i\rho,\qquad w=x^{3/2}=\rho^{3/2}e^{3\pi i/4},
 \qquad y=-w=\rho^{3/2}e^{-\pi i/4}.
\]
This point belongs to every such product sector and satisfies
$y^2=x^3$, but
\[
 f(x,y)=H(-w)-H(w)=-\frac\pi w e^{-i/w}\ne0.
\]
If a holomorphic extension of $g$ existed, the identity $f=Pg$
would continue from the positive quadrant to a common smaller
product sector and force $f(x,y)=0$ at this point, a contradiction.
\end{proof}

Thus exact division fails because the quotient loses the required
sectorial continuation. Estimate~\eqref{eq:D-real-bound} shows that
it retains the real Gevrey bounds: this example does \emph{not}
exhibit a loss of the numerical Gevrey order.

\end{document}